\documentclass[nonblindrev,opre]{informs3generic} 

\OneAndAHalfSpacedXI

\usepackage{endnotes}
\let\footnote=\endnote
\let\enotesize=\normalsize
\def\notesname{Endnotes}%
\def\makeenmark{$^{\theenmark}$}
\def\enoteformat{\rightskip0pt\leftskip0pt\parindent=1.75em
  \leavevmode\llap{\theenmark.\enskip}}

\usepackage[norelsize,linesnumbered,ruled,vlined]{algorithm2e}
\usepackage{booktabs}
\usepackage{multirow}
\usepackage{dsfont}
\usepackage{placeins}
\usepackage[caption=false,lofdepth,lotdepth]{subfig}
\usepackage{graphicx}

\makeatletter
\newcommand{\leqnomode}{\tagsleft@true\let\veqno\@@leqno}
\newcommand{\reqnomode}{\tagsleft@false\let\veqno\@@eqno}
\makeatother

\usepackage{natbib}
 \bibpunct[, ]{(}{)}{,}{a}{}{,}%
 \def\bibfont{\small}%
 \def\bibsep{\smallskipamount}%
 \def\bibhang{24pt}%
\TheoremsNumberedThrough     
\ECRepeatTheorems

\EquationsNumberedThrough    

\MANUSCRIPTNO{OPRE-2015-07-407.R2} 
                 
\let\chapter\section
 
\begin{document}


\RUNAUTHOR{Aswani, Shen and Siddiq}

\RUNTITLE{Inverse Optimization with Noisy Data}

\TITLE{Inverse Optimization with Noisy Data}

\ARTICLEAUTHORS{%
\AUTHOR{Anil Aswani, Zuo-Jun Max Shen, Auyon Siddiq}
\AFF{Department of Industrial Engineering and Operations Research, University of California, Berkeley, USA, \EMAIL{aaswani@berkeley.edu}, \EMAIL{maxshen@berkeley.edu}, \EMAIL{auyon.siddiq@berkeley.edu}.} 


} 

\ABSTRACT{%
Inverse optimization refers to the inference of unknown parameters of an optimization problem based on knowledge of its optimal solutions. This paper considers inverse optimization in the setting where measurements of the optimal solutions of a convex optimization problem are corrupted by noise. We first provide a formulation for inverse optimization and prove it to be NP-hard. In contrast to existing methods, we show that the parameter estimates produced by our formulation are statistically consistent. Our approach involves combining a new duality-based reformulation for bilevel programs with a regularization scheme that smooths discontinuities in the formulation. Using epi-convergence theory, we show the regularization parameter can be adjusted to approximate the original inverse optimization problem to arbitrary accuracy, which we use to prove our consistency results. Next, we propose two solution algorithms based on our duality-based formulation. The first is an enumeration algorithm that is applicable to settings where the dimensionality of the parameter space is modest, and the second is a semiparametric approach that combines nonparametric statistics with a modified version of our formulation. These numerical algorithms are shown to maintain the statistical consistency of the underlying formulation.  Lastly, using both synthetic and real data, we demonstrate that our approach performs competitively when compared with existing heuristics.
}%


\maketitle

%

\vspace{-10mm}
\section{Introduction}

An appreciable share of real-world data represents {\it decisions}, which can often be characterized as the solutions of correspondingly-defined optimization problems. Estimating the parameters of these latent optimization problems has the potential to provide greater insight into how decisions are made, and also enable the prediction of future decisions. Examples of domains where this is important include health systems engineering \citep{aswani2016}, energy systems engineering \citep{ratliff2014}, and marketing \citep{green1990conjoint}, where such estimation may lead to new approaches that enable the individualization of products and incentives. For example, consider a single homeowner who \emph{each day} observes an electricity price and weather forecast and then adjusts the temperature set-point for their home's air-conditioner.  By modeling this homeowner's decision as being generated from an optimization problem, we can directly estimate the price elasticity of comfort -- as measured by a standardized function of the temperature set-point and the outside temperature \citep{ashrae2013} -- for this particular homeowner.  This information is valuable for designing personalized incentive bonus schemes that encourage participation in demand-response programs \citep{aalami2010demand} or promote energy-efficiency \citep{aswani2012_allerton}.

\subsection{Overview}

This paper considers the problem of estimating unknown model parameters of an optimization problem based on noisy measurements of its optimal solutions, which is often referred to as {\it inverse} optimization.  In particular, we provide the first \emph{statistical inference} perspective on the inverse optimization problem.  This is important because real-world decision data is noisy, either because (i) the data collection process introduces measurement noise, (ii) the decision-maker deviates from optimal decisions -- phenomena often referred to as \emph{bounded rationality} \citep{tversky1981}, or (iii) there is mismatch between the parametric form of the model and the true underlying decision-making process.

Noisy data make inverse optimization challenging because noise in the solution data can preclude the existence of a single set of model parameters that renders all observed solutions exactly optimal. In this setting, the goal of inverse optimization is to find a set of model parameters that achieves a good ``fit" with respect to the solution data. More specifically, we are interested in two statistical questions.  First, how can we generate estimates of unknown model parameters that asymptotically provide the best possible predictions from the chosen parametric form of the model?  In statistics, this property is known as \emph{risk consistency} \citep{bartlett2002,greenshtein2004,chatterjee2014}.  Second, when the chosen model matches the true model that is is generating the solution data, how can we generate estimates that asymptotically converge to the true value of the unknown parameters?  In statistics, this property is known as \emph{consistency} \citep{wald1949,jennrich1969,bickel2006}.  We will use the term \emph{estimation consistency} to distinguish this concept from risk consistency.  Note that estimation consistency generally implies risk consistency.

Restated, a risk consistent estimate asymptotically achieves the lowest possible prediction error (out of all possible predictions permitted by the class of models considered).  Hence, risk consistency and estimation consistency allow us to be confident that prediction and estimation accuracy, respectively, will generally improve with additional data. By contrast, an estimator that fails to be risk consistent (so-called \emph{inconsistent} estimators) may yield poor predictions, even if a large amount of data is available.  Proving consistency of an estimator is an important topic in the theory of statistical inference (cf. \citep{wald1949,jennrich1969,bartlett2002,greenshtein2004,bickel2006,chatterjee2014,aswani2015}), and consistency is considered to be a minimal requirement for an estimator \citep{bickel2006}.

The main paper begins with Section \ref{section:challenges}, which describes the statistical and computational challenges of inverse optimization with noisy data.  The section begins by formally defining a (convex) forward optimization problem and its corresponding inverse optimization problem.  We specifically formulate the inverse optimization problem such that its solution has the desired statistical consistency properties. Our approach is conceptually similar to least squares regression in the sense that we also employ a sum-of-squares loss function to fit a parametric model to noisy data. The substantive difference is that inverse optimization involves estimating the (possibly multi-valued) solution set of a general convex optimization problem, whereas regression typically involves estimating a (single-valued) function which has a closed form expression. Due to these differences, much of the classical statistical theory on least-squares regression \citep{jennrich1969} is invalid in the inverse optimization setting, and thus new analysis is required. We also note that our approach is not restricted to the use of an $\ell_2$ norm: Results similar to those in our paper can be proved for other loss functions, such as absolute deviation or a likelihood function, but we do not consider those extensions in this paper.

In Section \ref{section:consistency}, we prove that our inverse optimization formulation produces  statistically consistent estimates of the unknown model parameters.  The key technical difficulty in proving these results is dealing with continuity issues.  In particular, the risk measures are not continuous in the general case, but are rather lower semicontinuous.  As alluded to above, this precludes the use of the typical statistical machinery used to prove consistency results (namely the uniform law of large numbers \citep{jennrich1969} and related uniform bounds \citep{bartlett2002,greenshtein2004}).  To circumvent this difficulty, we define a regularized version of the inverse optimization problem that smooths out any discontinuities, and this regularized version of the problem is constructed using a new duality-based reformulation for bilevel programs.  Using epi-convergence theory, we show the regularization parameter can be adjusted to approximate the original inverse optimization problem to arbitrary accuracy.  The regularized version of our formulation enables us to prove the desired statistical consistency results.

Section \ref{section:numerical} provides two numerical algorithms for solving our formulation of the inverse optimization problem.  The first numerical algorithm is an enumeration algorithm that is applicable to settings where the dimensionality of the parameter space is modest (i.e., at most four or five parameters).  The second numerical algorithm is a semiparametric approach that combines nonparametric statistics with a modified version of our formulation of the inverse optimization problem.  The statistical consistency of these two numerical algorithms are shown using the results from Section \ref{section:consistency}.  Lastly, in Section \ref{section:data} we demonstrate using synthetic and real data sets the competitiveness of our approaches as compared to existing heuristics \citep{keshavarz2011,bertsimas2013}.


\subsection{Literature Review}
\label{sec:litreview}
Existing inverse optimization models differ based on their specification of the \emph{loss function}, and the different models can be broadly categorized into either (i) deterministic settings, or (ii) noisy settings.  The work in the deterministic setting has primarily focused on single observation situations, wherein a single optimal solution is observed and then used to estimate parameters of the optimization problem.  However, in the noisy setting past work has considered situations with either a single observation and multiple observations. 

We begin by describing some of the work in the deterministic setting: \cite{ahuja2001} consider the estimation of objective function coefficients of general linear programs given a single optimal solution. The feasible region of the inverse problem is formulated using the constraints of the dual program and complemetary slackness conditions. Since the observed solution is assumed to be optimal, feasibility of the inverse problem is guaranteed. \cite{iyengar2005} and \cite{zhang2010} extend inverse optimization to certain conic forward problems using conic duality theory. Inverse optimization models have also been studied in the context of integer programs \citep{schaefer2009,wang2009} and network problems \citep{burton1992,hochbaum2003,zhang1996}. With respect to applications, inverse optimization models has been employed in many different domains, including healthcare \citep{erkin2010,chan2014}, energy \citep{ratliff2013incentive,saez2016}, finance \citep{bertsimas2012}, production planning \citep{troutt2006}, demand management \citep{carr2000,bajari2007}, auction design \citep{beil2003}, telecommunication \citep{farago2003} and geoscience \citep{burton1992}. We refer the reader to \cite{heuberger2004} for a survey of inverse optimization methods.

The noisy setting has been less studied.  \cite{chan2014} propose a generalized approach to inverse optimization for linear programs where the (single) observed solution may be suboptimal or infeasible.  Instead of complementary slackness, the authors use dual feasibility and strong duality to formulate the inverse problem. To accommodate noise, the strong duality constraint is relaxed to guarantee feasibility of the inverse problem.  \cite{saez2016} also consider inverse optimization for linear programs, and formulate the inverse problem using KKT conditions.  \cite{keshavarz2011} formulates the inverse problem using the KKT conditions of the optimization problem. To accommodate noise, the KKT conditions are relaxed by introducing slack variables to allow the data to  ``approximately'' satisfy the KKT conditions. Similarly, \cite{bertsimas2013} consider inverse problems where the observed data are assumed to be in an equilibrium. The authors enforce optimality conditions using a variational inequality, and similarly relax the optimality conditions by introducing slack variables to allow the data to ``approximately'' satisfy the variational inequality. 

Our work in this paper is most closely related to the noisy setting with multiple observations that has been previously considered by \cite{keshavarz2011} and \cite{bertsimas2013}. The key distinction between our work and these two previous approaches is in the choice of the loss function.  In \citep{keshavarz2011} and \citep{bertsimas2013}, the loss function is measured by the amount of slack required to make the measured data satisfy an approximate optimality condition (either the KKT conditions \citep{keshavarz2011} or a variational inequality describing optimality \citep{bertsimas2013}).  In contrast, our approach is to jointly estimate $(i)$ the parameters of the optimization problem, and $(ii)$ the denoised versions of the measured data (i.e. the true underlying optimal solutions). By performing this joint estimation, we are able to define our loss function to be the average discrepancy between the measured data and the (estimated) denoised data.  As we will show, this difference in loss function leads to significantly improved statistical performance.  A secondary distinction is that we propose the use of a novel optimality condition: specifically, we upper bound the objective function of a convex optimization problem by its dual -- thereby enforcing a zero duality gap and guaranteeing optimality.  An important benefit of using this alternate optimality condition is that it has favorable convexity and continuity properties (which are not available when using KKT conditions or variational inequalities to represent optimality) that enable design of numerical algorithms for solving the inverse optimization problem.

\subsection{Contributions}

Our contributions in this paper include both statistical and optimization results, and there are specifically two main contributions.  The first is we show that solving a bilevel formulation for the problem of inverse optimization with noisy data provides parameter estimates that are statistically consistent. This statistical result is independent of the approach used to solve the bilevel formulation.  Our second main contribution is to propose two numerical algorithms for solving the bilevel formulation by using a novel duality-based reformulation.  However, other numerical algorithms can be used to solve the bilevel formulation.  For instance, the bilevel program can be reformulated as a mixed-integer quadratic program (MIQP) in some cases \citep{fortuny1981,audet1997}.  Our statistical results apply to any numerical algorithm for solving the bilevel formulation, including the MIQP reformulation (when possible) or our two algorithms.

We also prove that existing heuristics for inverse optimization with noisy data \citep{keshavarz2011,bertsimas2013}, which are expressed as convex optimization problems, are statistically inconsistent -- meaning that in the limit of increasing amount of data these approaches will generate parameter estimates that converge to incorrect values. This is perhaps not unexpected, because we also prove that the problem of inverse optimization with noisy data is NP-hard.  It should be noted that the inverse optimization problem \emph{without} noisy data can be solved in polynomial time, as shown by \cite{keshavarz2011} and \cite{bertsimas2013}.

An additional contribution is we propose a novel reformulation of bilevel programs where there lower level optimization problem is convex.  It is common to replace the lower level problem by the KKT conditions or to upper bound the objective function by the value function \citep{dempe2015}.  However, these approaches face certain numerical difficulties.  We propose to upper bound the objective function by its dual, which enforces a zero duality gap and describes an optimal point.  The benefit of our optimality condition is it has convexity and continuity properties that support the design of numerical algorithms.  The two numerical algorithms we propose directly make use of this optimality condition, and the proofs of our statistical results are also aided by the use of this optimality condition.

\subsection{Notation}
Most notation we use in this paper is standard, and we briefly summarize some of the less usual aspects of our notation.  We use $\|\cdot\|$ to denote the usual $\ell_2$-norm.  The indicator function $\mathds{1}(p)$ is defined to be
\begin{equation}
\mathds{1}(p) = \begin{cases} 1, & \text{if condition } p \text{ is satisfied}\\
0, &\text{otherwise}\end{cases}
\end{equation}
The notation $[r] = \{1,\ldots,r\}$ refers to sequential set.  The Kuratowski limit superior of a sequence of sets $\mathcal{C}_\nu\subseteq\mathbb{R}^d$ is defined as
\begin{equation}
\textstyle\lim\sup_\nu (\mathcal{C}_\nu) = \{x \in\mathbb{R}^d : \liminf_\nu \text{dist}(x,\mathcal{C}_\nu) = 0\},
\end{equation}
where $\text{dist}(x,\mathcal{C}) = \inf\{\|x-c\|\ |\ c\in \mathcal{C}\}$.  We similarly define $\text{dist}(\mathcal{B},\mathcal{C}) = \inf\{\text{dist}(x,\mathcal{C})\ |\ x\in\mathcal{B}\}$. 

\section{Challenges with Noisy Inverse Optimization}
\label{section:challenges}

This section begins by formalizing the notation for the forward problem, before defining the noisy inverse optimization problem.  For the case where we have access to measurements (rather than the underlying distributions), we formulate a related sample average approximation of the inverse optimization problem.  We show that both these inverse problems are NP-hard.  We conclude by showing that existing heuristic approaches for solving the inverse optimization problem are statistically inconsistent, meaning that in the limit of infinite data these heuristic approaches converge to incorrect solutions.

\subsection{Model for Forward Problem}

Let $x\in\mathbb{R}^d$ be the decision variable, $u\in\mathbb{R}^m$ be the external input variable, and $\theta\in\mathbb{R}^p$ be the parameter vector.  Then the forward optimization problem is given by
\leqnomode\begin{equation*}
\tag*{$\mathsf{FOP}$}\min_x \big\{f(x, u,\theta)\ \big|\ g(x, u,\theta) \leq 0\big\},
\end{equation*}\reqnomode
where $f : \mathbb{R}^d \times \mathbb{R}^m \times \mathbb{R}^p \rightarrow \mathbb{R}$ is a function and $g : \mathbb{R}^d \times \mathbb{R}^m \times \mathbb{R}^p \rightarrow \mathbb{R}^q$ is a vector-valued function.
The solution set of $\mathsf{FOP}$ is the set-valued function given by $\mathcal{S}(u,\theta) = \arg\min_x \{f(x, u,\theta)\ |\ g(x, u,\theta) \leq 0\}$.
The value function of $\mathsf{FOP}$ is given by $V(u,\theta) = \min_x \{f(x,u,\theta)\ |\ g(x,u,\theta)\leq0\}$, and the feasible set is defined as $\Phi(u,\theta) = \big\{x \in\mathbb{R}^d : g(x,u,\theta) \leq 0\big\}$.


\subsection{Model for Inverse Optimization Problem}

Suppose $(u,y) \in \mathbb{R}^m\times\mathbb{R}^d$ is a vector-valued random variable that is distributed according to some unknown but fixed joint distribution $\mathbb{P}_{(u,y)}$.  Let $\mathcal{U}\times\mathcal{Y}\subseteq\mathbb{R}^m$ be the support of this distribution, meaning the smallest set that satisfies the property $\mathbb{P}_{(u,y)}(\mathcal{U},\mathcal{Y}) = 1$.  If we define the function
\leqnomode\begin{equation}
\tag*{$\mathsf{RISK}$}Q(\theta) = \mathbb{E}\Big(\min_{x\in\mathcal{S}(u,\theta)}\|y - x\|^2\Big),
\end{equation}\reqnomode
then the inverse optimization problem is given by
\leqnomode\begin{equation*}
\tag*{$\mathsf{IOP}$}\min \big\{Q(\theta)\ \big|\ \theta\in\Theta\big\},
\end{equation*}\reqnomode
where $\Theta \subseteq\mathbb{R}^p$ is a known set.  We make the following assumptions:\\

\noindent\textbf{A1. } The functions $f(x, u,\theta)$ and $g(x, u,\theta)$ are continuous in $x,u,\theta$ and convex in $x$ for fixed $u,\theta$.\\

\noindent\textbf{A2. } The set $\Theta$ is convex.\\

\noindent These assumptions are fairly mild. \textbf{A1} is equivalent to stating $\mathsf{FOP}$ is a convex optimization problem. Though \textbf{A2} is necessary for the semiparametric algorithm presented in Section 4 because it ensures polynomial time computability of the algorithm, it is not necessary for our main results regarding statistical consistency because these results only require that $\Theta$ is well-posed.  Hence, \textbf{A2} is one way to ensure $\Theta$ is well-posed, and one alternative for which our statistical consistency results would hold is if $\Theta$ is discrete-valued and finite.

When the joint distribution $\mathbb{P}_{(u,y)}$ is unknown, we cannot solve $\mathsf{IOP}$ without additional information.  Fortunately, we can leverage the iid measurements $(u_i,y_i)$ for $i \in [n]$. In principle, we can solve $\mathsf{IOP}$ using a sample average approximation:
\leqnomode\begin{equation}
\tag*{$\mathsf{IOP}$--$\mathsf{SAA}$}\min\big\{Q_n(\theta)\ \big|\ \theta\in\Theta\big\},
\end{equation}\reqnomode
where
\leqnomode\begin{equation}
\tag*{$\mathsf{RISK}$--$\mathsf{SAA}$}
\begin{aligned}
Q_n(\theta) = \min_{x_i}\ &\frac{1}{n}\sum_{i=1}^n\|y_i-x_i\|^2\\
\text{s.t. }&x_i \in \mathcal{S}(u_i,\theta),&\forall i\in[n]
\end{aligned}
\end{equation}\reqnomode
In the context of a decision-making agent, $u_i$ may be interpreted as an external signal the agent responds to and $y_i$ as a noisy observation of the corresponding decision of the agent. Note that in the expression \textsf{RISK}, the variable $x$ is constrained to be an optimal solution of the forward problem. Similarly, we may interpret $x_i$ as representing  an underlying optimal solution (unperturbed by noise) of \textsf{FOP} in the $i^{th}$ instance. Note also that while the $u_i$ and $\theta$ are both parameters of \textsf{FOP}, they are different in that the $u_i$ are known and may vary across the $n$ observations, whereas $\theta$ is unknown and is fixed across all instances.

For a concrete example, consider the numerical experiments presented in Section 5.4, where we estimate an individual's utility function capturing the tradeoff between maintaing a comfortable indoor temperature versus the amount of energy consumption (and implicitly the air conditioning energy costs) required to cool the room. In that example, the $u$ represents the outside air temperature, $\theta_1$ captures the decision-maker's (unknown) tradeoff between comfort and energy consumption, $\theta_2$ parameterizes their (unknown) preferred temperature (i.e., the preferred temperature is $\theta_2 + u$), $x$ represents the true optimal temperature setpoint (for the given $u$ and $\theta$), and $y$ represents the temperature set-point that we observe.

\subsection{NP-Hardness of Inverse Optimization Problem}

Though all the functions and sets involved in $\mathsf{FOP}$ and $\mathsf{IOP}$ are convex, solving $\mathsf{IOP}$ is NP-hard.

\begin{theorem}
\label{thm:nphard}
If $\mathbf{A1}$,$\mathbf{A2}$ hold, then $\mathsf{IOP}$ is NP-hard.
\end{theorem}

\proof{Proof. }We prove this by showing a reduction from the problem of computing the best rank-1 approximation of an order 3 tensor (which is NP-hard \citep{hillar2013}) to $\mathsf{IOP}$.  Consider any $\psi \in \mathbb{R}^{r_1\times r_2\times r_3}$, where $r_1,r_2,r_3 \in \mathbb{Z}_+$.  This defines $\psi$ to be an order 3 tensor.  We define $\rho = r_1 + r_2 + r_3$, and suppose the parameter vector is given by $\theta = (a, b, c) \in \Theta = \mathbb{R}^\rho$, where $a\in\mathbb{R}^{r_1}$, $b\in\mathbb{R}^{r_2}$, and $c\in\mathbb{R}^{r_3}$.  Also define the discrete set $\mathcal{U} = [r_1]\times[r_2]\times[r_3]$, and suppose that $u = (\alpha,\beta,\gamma)$ is uniformly distributed over $\mathcal{U}$.  Furthermore, suppose $y$ is a random variable given by $\psi_{\alpha,\beta,\gamma}$, which means that $y$ is dependent on $u$ since $u = (\alpha,\beta,\gamma)$.  Then we define the following forward optimization problem
\begin{equation}
\mathcal{S}(u,\theta) = \arg\min_x \Big(x -  a_\alpha\cdot b_\beta\cdot c_\gamma\Big)^2.
\end{equation}
This forward optimization problem is a quadratic program (QP) when $(u,\theta)$ is fixed, and so the solution set is $\mathcal{S}(u,\theta) =  a_\alpha b_\beta c_\gamma$.  Note that the solution set consists of a single point.  Next, observe that
\begin{equation}
\min_{\theta\in\Theta}\ Q(\theta) = \min_{\theta\in\mathbb{R}^\rho} \frac{1}{\rho}\sum_{\alpha=1}^{r_1}\sum_{\beta=1}^{r_2}\sum_{\gamma=1}^{r_3}\Big(\psi_{\alpha,\beta,\gamma}- a_\alpha\cdot b_\beta\cdot c_\gamma\Big)^2,\label{eqn:rank1}
\end{equation}
where we have converted the expectation into a weighted sum using the fact that $u$ is uniformly distributed over $\mathcal{U}$.  Observe that (\ref{eqn:rank1}) is the problem of computing the best rank-1 approximation to an order 3 tensor \citep{hillar2013}.
\Halmos\endproof


\begin{remark}
\label{remark:iars}
Inapproximability results for $\mathsf{IOP}$ can be shown under the setting where $\Theta$ is allowed to be a discrete set (i.e, \textbf{A1} holds, but \textbf{A2} does not hold).  In particular, there is a straightforward reduction from the shortest vector problem.  This implies that $\mathsf{IOP}$ is NP-hard to approximate to within any factor up to $2^{{(\log d)}^{1-\epsilon}}$, for any $\epsilon \geq 0$ \citep{haviv2012}.
\end{remark}

\begin{remark}
\label{remark:saapts}
Polynomial-time solvability of $\mathsf{IOP}$ is possible in very specific settings.  For instance, if $\mathsf{FOP}$ is a QP with the solution set $\mathcal{S}(u,\theta) = \arg\min_x \{x^2 - 2(\theta + u)\cdot x\} = \theta + u$ or an LP with the solution set $\mathcal{S}(u,\theta) = \arg\min_x \{x : x = \theta + u\} = \theta + u$, then $\mathsf{IOP}$ is a QP: $\min_{\theta \in \Theta} \big\{\mathbb{E}((y - \theta - u)^2)\big\}$, and its minimizer is $\theta^* = \mathbb{E}(y-u)$.

\end{remark}
In general, since $\mathcal{S}(u_i,\theta)$ is the optimal solution sets to \textsf{FOP} under input $u_i$, the problem $\mathsf{IOP}$--$\mathsf{SAA}$ is a bilevel program, which are usually difficult to solve \citep{dempe2015}. In fact, $\mathsf{IOP}$--$\mathsf{SAA}$ is also NP-hard to solve.

\begin{remark}
In the case where \textsf{FOP} is a linear program, the inverse problem \textsf{IOP} takes the form of a quadratic bilevel program, which are generally NP-hard \citep{audet1997}. Branch-and-bound algorithms have been proposed for solving such bilevel programs \citep{bard1990branch}.
\end{remark}

\begin{corollary}
If $\mathbf{A1}$,$\mathbf{A2}$ hold, then $\mathsf{IOP}$--$\mathsf{SAA}$ is NP-hard.
\end{corollary}

\proof{Proof. }We show this result using the same construction used to prove Theorem \ref{thm:nphard}.  In particular, observe that if $\{u_1,\ldots,u_n\} = \mathcal{U}$, then $\mathsf{IOP}$--$\mathsf{SAA}$ is equivalent to $\mathsf{IOP}$ , which is NP-hard by Theorem \ref{thm:nphard}.  Finally, note that the condition $\{u_1,\ldots,u_n\} = \mathcal{U}$ occurs with nonzero probability since the set $\mathcal{U}$ is finite and since the $u_i$ are sampled uniformly from $\mathcal{U}$.
\Halmos\endproof

\begin{remark}
Inapproximability results for $\mathsf{IOP}$--$\mathsf{SAA}$ can be shown under the setting where $\Theta$ is allowed to be a discrete set (i.e, \textbf{A1} holds, but \textbf{A2} does not hold).  In particular, the same construction in Remark \ref{remark:iars} can be used to shown $\mathsf{IOP}$--$\mathsf{SAA}$ is NP-hard to approximate to within any factor up to $2^{{(\log d)}^{1-\epsilon}}$, for any $\epsilon \geq 0$ \citep{haviv2012}.
\end{remark}

\begin{remark}
Polynomial-time solvability of $\mathsf{IOP}$--$\mathsf{SAA}$ is possible in very specific settings.  For instance, the constructions in Remark \ref{remark:saapts} lead to instances of $\mathsf{IOP}$--$\mathsf{SAA}$ that are QP's.
\end{remark}

\subsection{Statistical Consistency in Inverse Optimization with Noisy Data}

We begin with two statistical definitions of consistency: risk consistency and estimation consistency.  These definitions are stated in order of increasing stringency, meaning that risk consistency is necessary (in situations with sufficient continuity) for estimation consistency.  The first definition relates to the best predictions possible using the given forward optimization problem.

\begin{definition}[Risk Consistency]
An estimate $\hat{\theta}_n \in \Theta$ is risk consistent if
\begin{equation}
Q(\hat{\theta}_n) \stackrel{p}{\longrightarrow} \min \big\{Q(\theta)\ \big|\ \theta\in\Theta\big\}.
\end{equation}
\end{definition}

\noindent We should interpret the function $Q(\theta)$ as the expected prediction error when the parameter values are $\theta$, where the prediction is the solution set $\mathcal{S}(u,\theta)$.  And so the above definition is stating that an estimator $\theta_n$ is risk consistent if the expected prediction error of the estimate $\theta_n$ converges in probability to the minimum prediction error possible when we use the forward optimization model described by \textsf{FOP} and constrain $\theta$ to belong to $\Theta$.  In other words, an estimator is risk consistent if it asymptotically provides the best predictions possible.

The second statistical definition relates to the situation where the forward optimization model described by \textsf{FOP} is correct and there is a \emph{true} parameter.  In particular, it applies to situations where the below identifiability condition is satisfied.  Briefly summarized, the identifiability condition is satisfied when \textsf{FOP} is such that two different parameter values $\theta_1$ and $\theta_2$ lead to two different distributions for measurements of the decision data $y_i$.  More details and clarifying examples are found in Appendix \ref{sect:iden}.\\

\noindent\textbf{IC. } There exists a unique $\theta_0\in \Theta$ such that the following three sub-conditions hold: (i) $y = \xi + w$, where $\xi \in \mathcal{S}(u,\theta_0)$, $\mathbb{E}(w) = 0$, $\mathbb{E}(w^2) < +\infty$, and $u,\xi$ are independent of $w$, (ii) for all $\theta \in \Theta\setminus\theta_0$ there exists $\mathcal{U}(\theta)\subseteq\mathcal{U}$ such that $\mathbb{P}(u \in \mathcal{U}(\theta)) > 0$ and $\text{dist}(\mathcal{S}(u,\theta), \mathcal{S}(u,\theta_0)) > 0$ for each $u\in\mathcal{U}(\theta)$, and (iii) for each fixed $\theta\in\Theta$ we have $\mathbb{P}(\{u : \mathcal{S}(u,\theta) \text{ is multivalued}\}) = 0$.\\


\noindent The first sub-condition of the identifiability condition is stating that the solution data $y_i$ is a noisy measurement (with noise random variable $w$) of a point that belongs to the solution set $\mathcal{S}(u_i,\theta_0)$, and the second sub-condition is stating that when $\theta$ is different from $\theta_0$ then this leads to different solution sets.  This second sub-condition is necessary, because otherwise we could not distinguish the predictions of \textsf{FOP} when the parameters $\theta$ differ from $\theta_0$.  The third sub-condition eliminates pathological cases that occur when the solution set at a fixed $\theta$ is so large that it approximately encompasses all possible solutions.  Note that this third sub-condition is mild, and examples where it is satisfied include when (i) \textsf{FOP} is strictly convex, or when (ii) \textsf{FOP} is a linear program with random coefficients drawn from a continuous distribution; it holds for other examples as well.  The second statistical definition is related to this identifiability condition.

\begin{definition}[Estimation Consistency]
Suppose $\mathbf{IC}$ holds.  An estimate $\hat{\theta}_n \in \Theta$ is estimation consistent if
\begin{equation}
\hat{\theta}_n \stackrel{p}{\longrightarrow} \theta_0.
\end{equation}
\end{definition}

\noindent Stated in words, an estimate $\hat{\theta}_n$ is estimation consistent if it converges in probability to the true parameter values $\theta_0$.  This is the classical notion of consistency of a statistical estimator \citep{bickel2006}.  

Though these statistical notions of consistency are quite natural, it is the case that existing heuristic approaches for solving the inverse optimization problem are statistically inconsistent.  We will use $\mathsf{VIA}$ to refer to the variational inequality method of \cite{bertsimas2013}, and we refer to the KKT conditions approach of \cite{keshavarz2011} as $\mathsf{KKA}$.

\begin{proposition}
\label{proposition:estincon}
Suppose $\mathbf{A1}$,$\mathbf{A2}$ and $\mathbf{IC}$ hold.  Then $\mathsf{VIA}$ \citep{bertsimas2013} and $\mathsf{KKA}$ \citep{keshavarz2011} are not estimation consistent.
\end{proposition}

\begin{corollary}
Suppose $\mathbf{A1}$,$\mathbf{A2}$ hold.  Then $\mathsf{VIA}$ \citep{bertsimas2013} and $\mathsf{KKA}$ \citep{keshavarz2011} are not risk consistent.
\end{corollary}

The proofs for Proposition 1 and Corollary 2 are contained in the Appendix. The intuition for why $\mathsf{VIA}$ and $\mathsf{KKA}$ are statistically inconsistent is that they are minimizing an incorrect measure of error: These approaches generate an estimated set of parameters that minimizes the level of suboptimality of the measured solution data.  However, this leads to biased estimates because suboptimality is measured by (i) deviations in the value of the objective function of \textsf{FOP} and (ii) the amount of constraint violation of \textsf{FOP}, whereas noise directly perturbs the solution data. This is in contrast to our approach (as exemplified by \textsf{IOP}-\textsf{SAA}) which generate an estimated set of parameters that minimizes the deviation between predicted and measured solution data. This distinction between suboptimality and deviations in the solution data becomes most apparent (and critical) in problems with constraints.

\section{Consistent Estimation for Inverse Optimization Problem}
\label{section:consistency}

Given the statistical inconsistency of existing heuristics, we propose to solve the noisy inverse optimization problem by instead solving \textsf{SAA}-\textsf{IOP}.  First, we will need to impose a regularity condition to ensure that $\mathsf{FOP}$ and $\mathsf{IOP}$--$\mathsf{SAA}$ are numerically well-posed:\\

\noindent\textbf{R1. } For each $u \in\mathcal{U}$ and $\theta\in\Theta$, the feasible set $\Phi(u,\theta)$ is closed, bounded, and has a nonempty interior (i.e., $\text{int}(\Phi(u,\theta)) \neq \emptyset$).  The feasible set $\Phi(u,\theta)$ is also absolutely bounded, meaning there exists $M > 0$ such that $\|x\|\leq M$, for all $x\in\Phi(u,\theta)$, $u\in\mathcal{U}$, and $\theta\in\Theta$.\\

\noindent Condition \textbf{R1} is equivalent to requiring $\mathsf{FOP}$ to have a strictly feasible point (i.e., Slater's condition holds), and that the feasible set of $\mathsf{FOP}$ is closed and bounded. The first sub-condition requiring the feasible set be closed and bounded is needed to ensure the existence of well-posed primal and dual solutions, and it could be replaced by more general conditions.  For instance, we could have instead assumed \textsf{FOP} satisfies the uniform level-boundedness condition \citep{rockafellar1998}. We use the above for simplicity of stating the results.  The condition that $\Phi(u,\theta)$ has a nonempty interior\footnote{}\footnotetext[1]{\textbf{R1} can be relaxed to requiring a nonempty \emph{relative} interior if the affine constraints of \textsf{FOP} are of the form $Mx + \zeta(u,\theta) = 0$, where $M$ is a matrix and $\zeta$ is a continuous function. The reason is that our proofs make use of a result (Example 5.10 of \cite{rockafellar1998}) on the continuity of parametrized convex constraints with a nonempty interior, and this result can be generalized for the above case through minor modifications (using corresponding results on relative interiors from Section 2.H of \cite{rockafellar1998}) to ensure continuity of the feasible set of \textsf{FOP} with a nonempty relative interior.  Generalizing Example 5.10 of \cite{rockafellar1998} or our results to cases with more complex affine constraints will require further study.} is needed to ensure continuity of $\mathcal{S}(u,\theta)$ through application of the Berge Maximum Theorem \citep{berge1963}.


The simplest case of statistical consistency of \textsf{SAA}-\textsf{IOP} occurs when the function $f(x,u,\theta)$ is strictly convex, because of the following result:

\begin{proposition}
\label{proposition:sccont}
Suppose $\mathbf{A1}$,$\mathbf{A2}$ and $\mathbf{R1}$ hold.  If $f(x,u,\theta)$ is strictly convex in $x$ for fixed $u\in\mathcal{U}$ and $\theta\in\Theta$, then $Q_n(\theta)$ is continuous.
\end{proposition}

\proof{Proof.} Because the feasible set $\Phi(u,\theta)$ is convex for fixed $u,\theta$ by \textbf{A1} and has a nonempty interior by \textbf{R1}, this means $\Phi(u,\theta)$ is continuous in $\theta$ by Example 5.10 from \citep{rockafellar1998}.  Thus, we can apply the Berge Maximum Theorem \citep{berge1963} to $\mathsf{FOP}$.  This implies $\mathcal{S}(u,\theta)$ is upper hemicontinuous in $\theta$ for fixed $u\in\mathcal{U}$.  However, $\mathcal{S}(u,\theta)$ consists of a single point for fixed $u\in\mathcal{U}$ and $\theta\in\Theta$, because the objective function is strictly convex and since $\mathbf{R1}$ holds.  Consequently, $\mathcal{S}(u,\theta)$ is a continuous single-valued function for fixed $u\in\Theta$ (see for instance Theorem 2.6 in \citep{rockafellar1998}).  Thus, we can apply the Berge Maximum Theorem to $\mathsf{RISK}$--$\mathsf{SAA}$, and this implies that $Q_n(\theta)$ as defined in $\mathsf{RISK}$--$\mathsf{SAA}$ is continuous.
\Halmos\endproof

In this case, we can prove risk and estimation consistency using standard arguments \citep{jennrich1969,van2000,bickel2006} from statistics that use the uniform law of large numbers \citep{jennrich1969}.  However, this approach cannot be applied to the more general case where $f(x,u,\theta)$ is not strictly convex.  In particular, when $f(x,u,\theta)$ is not strictly convex, the function $Q_n(\theta)$ will not generally be continuous.  And so a different argument is required because the uniform law of large numbers does not apply to discontinuous functions.  

Our approach will be to use a statistical consistency result originally due to \cite{wald1949} that uses a one-sided bounding argument.  The advantage of this approach is that it only requires lower semicontinuity, which we show always holds for $Q_n(\theta)$.  However, this result only implies the estimates $\hat{\theta}_n$ converge in probability to the set of minimizers of $Q(\theta)$.  This cannot imply risk consistency in the general case because $Q_n(\theta)$ is lower semicontinuous, which means that $Q(\hat{\theta}_n)$ can remain bounded from the minimum $Q(\theta)$.  And so for the general case, we will show that a weak risk consistency result holds.

To develop the statistical consistency results for the most general case, we will develop a regularized version of $\mathsf{RISK}$--$\mathsf{SAA}$ that is guaranteed to be continuous.  The first step of this construction involves proposing a new reformulation for bilevel programs that we call a duality-based reformulation.  Next, we use this reformulation to construct a regularized version of $\mathsf{RISK}$--$\mathsf{SAA}$ and prove its continuity.  We use this regularized version to prove statistical consistency results about $\mathsf{IOP}$--$\mathsf{SAA}$ and a regularized version of $\mathsf{IOP}$--$\mathsf{SAA}$.

\subsection{Duality-Based Reformulation}

One approach to solving bilevel problems (such as $\mathsf{IOP}$--$\mathsf{SAA}$) is to reformulate the problem as a normal (i.e., single level) optimization problem by replacing the constraints $x_i \in \mathcal{S}(u_i,\theta)$ with an optimality condition \citep{dempe2015}.  One possibility is to replace $x_i \in \mathcal{S}(u_i,\theta)$ by the KKT conditions of $\mathsf{FOP}$, and another possibility is to upper bound the objective function using the value function $f(x_i,u_i,\theta) \leq V(u_i,\theta)$.  Unfortunately, these approaches often encounter numerical difficulties.  The KKT approach leads to a nonlinear program with combinatorial complexity, because of the complimentary slackness in KKT. The value function approach is difficult to implement because closed-form expressions for the value function are not available except for very special cases.

Here, we present a new optimality condition.  Given the numerical difficulties of existing approaches, we propose to solve bilevel programs (such as $\mathsf{IOP}$--$\mathsf{SAA}$) by using the Lagrangian dual function to upper bound the objective function.  The following proposition shows that our idea of using the dual as an upper bound represents a novel optimality condition.

\begin{proposition}
\label{prop:optref}
Suppose $\mathbf{A1}$ and $\mathbf{R1}$ hold.  Then $x \in \mathcal{S}(u,\theta)$ if and only if there exists a corresponding $\lambda \in \mathbb{R}^q$ for which $x,\lambda$ satisfy the inequalities
\begin{equation}
\label{eqn:dual_ineq}
\begin{aligned}
&f(x,u,\theta)-h(\lambda,u,\theta)\leq 0\\
&g(x,u,\theta)\leq0\\
&\lambda\geq0
\end{aligned}
\end{equation}
where $h(\lambda,u,\theta)$ is the Lagrangian dual function of $\mathsf{FOP}$. 
\end{proposition}



Proposition 3 is a consequence of strong duality for convex optimization problems. We can now exactly reformulate $\mathsf{RISK}$--$\mathsf{SAA}$ as the following optimization problem:
\leqnomode\begin{equation}
\tag*{$\mathsf{DB}$--$\mathsf{RISK}$--$\mathsf{SAA}$}
\begin{aligned}
Q_n(\theta) = \min_{x_i,\lambda_i}\ &\frac{1}{n}\sum_{i=1}^n\|y_i-x_i\|^2\\
\text{s.t. }&f(x_i,u_i,\theta)-h(\lambda_i,u_i,\theta)\leq 0,&\forall i\in[n]\\
&g(x_i,u_i,\theta)\leq0,&\forall i\in[n]\\
&\lambda_i\geq0,&\forall i\in[n]
\end{aligned}
\end{equation}\reqnomode
It should be noted that the formulation $\mathsf{DB}$--$\mathsf{RISK}$--$\mathsf{SAA}$ requires the Lagrangian dual function $h(\lambda,u,\theta)$ to be computable in closed form, which is the case for a large class of convex (e.g., linear, quadratic, conic) optimization problems that arise in practice \citep{boyd2009}. In cases where the dual function does not have an analytical representation, we may still solve $\mathsf{DB}$--$\mathsf{RISK}$--$\mathsf{SAA}$ by developing an algorithm that computes $h(\lambda,u,\theta)$ numerically, although designing such an algorithm is beyond the scope of this paper.  

One important feature of this reformulation is that it is a convex optimization problem for fixed values of $\theta$.

\begin{proposition}
\label{prop:convex}
Suppose $\mathbf{A1}$ and $\mathbf{R1}$ hold.  Then $\mathsf{DB}$--$\mathsf{RISK}$--$\mathsf{SAA}$ is a convex optimization problem for fixed $\theta$. 
\end{proposition}
Proposition 4 follows directly from ${\bf A1}$ and the concavity of the dual function in $\lambda$.


\subsection{Regularized Formulation}

Recall that $Q_n(\cdot)$ is generally not continuous even when $\mathbf{A1}$,$\mathbf{A2}$,$\mathbf{R1}$ hold.  Consequently, we develop a regularized version of the duality-based problem that is guaranteed to be continuous.  We define the $\epsilon$-regularized version of the duality-based problem to be
\leqnomode\begin{equation}
\tag*{$\mathsf{R}$--$\mathsf{DB}$--$\mathsf{RISK}$--$\mathsf{SAA}$}
\begin{aligned}
Q_n(\theta; \epsilon) = \min_{x_i,\lambda_i}\ &\frac{1}{n}\sum_{i=1}^n\|y_i-x_i\|^2\\
\text{s.t. }&f(x_i,u_i,\theta)-h(\lambda_i,u_i,\theta)\leq \epsilon,&\forall i\in[n]\\
&g(x_i,u_i,\theta)\leq\epsilon,&\forall i\in[n]\\
&\lambda_i\geq0,&\forall i\in[n]
\end{aligned}
\end{equation}\reqnomode
We associate this to a regularized version of the sample average approximation of the inverse optimization problem:
\leqnomode\begin{equation}
\tag*{$\mathsf{R}$--$\mathsf{IOP}$--$\mathsf{SAA}$}
\min\{Q_n(\theta; \epsilon)\ |\ \theta\in\Theta\}.
\end{equation}\reqnomode
The idea of this regularization is that we relax the optimality conditions to allow points $x_i$ to be an $\epsilon$-optimal solution.  Recall that a point
\begin{equation}
x^\epsilon \in \epsilon\text{-}\arg\min\{f(x)\ |\ g(x) \leq 0\},
\end{equation}
if (i) $f(x^\epsilon)-f^* \leq \epsilon$ and (ii) $g(x^\epsilon) \leq \epsilon$, where $f^*= \min\{f(x)\ |\ g(x) \leq 0\}$.

\begin{proposition}
Suppose $\mathbf{A1}$ and $\mathbf{R1}$ hold.  Then a point $x$ is an $\epsilon$-optimal solution if and only if there exists a corresponding $\lambda \in \mathbb{R}^q$ for which $x,\lambda$ satisfy the inequalities
\begin{equation}
\label{eqn:dual_ineq_e}
\begin{aligned}
&f(x,u,\theta)-h(\lambda,u,\theta)\leq \epsilon\\
&g(x,u,\theta)\leq\epsilon\\
&\lambda\geq0
\end{aligned}
\end{equation}
where $h(\lambda,u,\theta)$ is the Lagrangian dual function of $\mathsf{FOP}$. 
\end{proposition}



One benefit of this regularization is that it ensures convexity of $\mathsf{R}$--$\mathsf{DB}$--$\mathsf{RISK}$--$\mathsf{SAA}$ when $\theta$ is fixed.

\begin{proposition}
\label{prop:convex_r}
Suppose $\mathbf{A1}$,$\mathbf{A2}$ and $\mathbf{R1}$ hold.  Then $\mathsf{R}$--$\mathsf{DB}$--$\mathsf{RISK}$--$\mathsf{SAA}$ is a convex optimization problem for fixed $\theta$. 
\end{proposition}


Though the above propositions show that the regularization is equivalent to replacing optimality conditions with $\epsilon$-optimality conditions while maintaining convexity for fixed values of $\theta$, the main benefit of the regularization is that it ensures the function $Q_n(\theta; \epsilon)$ defined in $\mathsf{R}$--$\mathsf{DB}$--$\mathsf{RISK}$--$\mathsf{SAA}$ is continuous in $\theta,\epsilon$ for any $\epsilon > 0$.

\begin{proposition}
\label{prop:cont}
Suppose $\mathbf{A1}$,$\mathbf{A2}$ and $\mathbf{R1}$ hold.  Then the function $Q_{n}(\theta; \epsilon)$ is jointly continuous in $\theta, \epsilon$ for any $\epsilon > 0$.
\end{proposition}

\proof{Proof. }The solution set $\mathcal{S}(u,\theta)$ is nonempty under $\mathbf{A1}$,$\mathbf{R1}$ (see for instance Theorem 1.9 of \citep{rockafellar1998}).  Pick any $x_i \in \mathcal{S}(u_i,\theta)$, and let $\lambda_i$ be such that $x_i,\lambda_i$ satisfy (\ref{eqn:dual_ineq}) -- this $\lambda_i$ exists by Proposition \ref{prop:optref}.  Next, consider the sets
\begin{equation}
\begin{aligned}
&\overline{\mathcal{S}}(u_i,\theta; \epsilon) = \{x : f(x,u_i,\theta)-h(\lambda_i,u_i,\theta) \leq \epsilon,\ g(x,u_i,\theta) \leq \epsilon\}\\
&\mathcal{S}(u_i,\theta; \epsilon) = \{x : f(x,u_i,\theta)-h(\lambda,u_i,\theta) \leq \epsilon,\ g(x,u_i,\theta) \leq \epsilon,\ \lambda \geq 0\}
\end{aligned}
\end{equation}
and note that $\overline{\mathcal{S}}(u_i,\theta; \epsilon) = \mathcal{S}(u_i,\theta; \epsilon)$, since by optimality of $\lambda_i$ with respect to the dual problem we have $h(\lambda,u_i,\theta) \leq h(\lambda_i,u_i,\theta)$ for all $\lambda \geq 0$.  Observe that the functions $f(x_i,u_i,\theta), g(x_i,u_i,\theta)$ are continuous and convex from $\mathbf{A1}$, and the point $x_i$ belongs to the interior of $\overline{\mathcal{S}}(u_i,\theta; \epsilon)$ since it satisfies (\ref{eqn:dual_ineq}).  Thus, we can apply Example 5.10 from \citet{rockafellar1998}: This yields that $\overline{\mathcal{S}}(u_i,\theta; \epsilon)$ is continuous in $\theta,\epsilon$ for any $\epsilon > 0$, and so we also get continuity of $\mathcal{S}(u_i,\theta; \epsilon)$ by its equality to $\overline{\mathcal{S}}(u_i,\theta; \epsilon)$.  Since $\mathsf{R}$--$\mathsf{DB}$--$\mathsf{RISK}$--$\mathsf{SAA}$ can be written as $Q_n(\theta; \epsilon) = \min_{x_i}\{\frac{1}{n}\sum_{i=1}^n\|y_i-x_i\|^2\ |\ x_i \in \mathcal{S}(u_i,\theta; \epsilon),\ \forall i\in[n]\}$, we are able to apply the Berge Maximum Theorem \citep{berge1963}.  This implies continuity of $Q_n(\theta; \epsilon)$ in $\theta,\epsilon$ for any $\epsilon > 0$.
\Halmos\endproof

A point of note is that within the above proof, we show that the set of $\epsilon$-optimal solutions of a parametric convex optimization problem $\mathcal{S}(u_i,\theta; \epsilon)$ is continuous with respect to the parametrization $\theta$; this is in contrast to the solution set of a parametric convex optimization problem $\mathcal{S}(u_i,\theta)$, which is in general only upper hemicontinuous with respect to the parametrization $\theta$.  The case of a parametric strictly convex optimization problem is the exception, which as shown in the proof of Proposition \ref{proposition:sccont} has a continuous (with respect to the parametrization $\theta$) solution set.

The function $Q_n(\theta; \epsilon)$ will not be jointly continuous in $\theta, \epsilon$ at $\epsilon = 0$.  However, it satisfies another property that is useful for solving $\mathsf{IOP}$--$\mathsf{SAA}$:

\begin{proposition}
\label{prop:convmin}
Suppose $\mathbf{A1}$,$\mathbf{A2}$ and $\mathbf{R1}$ hold, and let $\epsilon_\nu > 0$ be a monotone decreasing sequence with $\epsilon_\nu\rightarrow 0$.  Then we have $\min \{Q_n(\theta; \epsilon_\nu)\ |\ \theta\in\Theta\} \rightarrow\min \{Q_n(\theta)\ |\ \theta\in\Theta\}$ and 
\begin{equation}
\textstyle\lim\sup_\nu (\arg\min \{Q_n(\theta; \epsilon_\nu)\ |\ \theta\in\Theta\}) \subseteq \arg\min \{Q_n(\theta)\ |\ \theta\in\Theta\}.
\end{equation}
If $z_\nu > 0$ is a monotone decreasing sequence with $z_\nu\rightarrow 0$, then we also have
\begin{equation}
\textstyle\lim\sup_\nu (z_\nu\text{-}\arg\min \{Q_n(\theta; \epsilon_\nu)\ |\ \theta\in\Theta\}) \subseteq \arg\min \{Q_n(\theta)\ |\ \theta\in\Theta\}.
\end{equation}
\end{proposition}

\proof{Proof. }Let $\mathcal{C}_n(\theta, \epsilon)$ be the feasible set of $\mathsf{R}$--$\mathsf{DB}$--$\mathsf{RISK}$--$\mathsf{SAA}$, and define $(X,\Lambda) = \{x_i,\lambda_i, \forall i\in[n]\}$.  Suppose $(X,\Lambda) \in \mathcal{C}_n(\theta, \alpha)$, where $\alpha \geq 0$.  Then for any $\beta \geq \alpha$ we must have $(X,\Lambda) \in \mathcal{C}_n(\theta, \beta)$ by the definition of the constraints in $\mathsf{R}$--$\mathsf{DB}$--$\mathsf{RISK}$--$\mathsf{SAA}$.  This means that
\begin{equation}
\mathcal{C}_n(\theta, \epsilon_1) \supseteq \mathcal{C}_n(\theta, \epsilon_2) \supseteq \cdots
\end{equation}
As a result, the set $\mathcal{D}_n(\theta, \epsilon_\nu) = \{\theta, X,\Lambda : \theta\in\Theta \text{ and } (X,\Lambda) \in \mathcal{C}_n(\theta,\epsilon_\nu)\}$ is also monotone nonincreasing:
\begin{equation}
\mathcal{D}_n(\theta, \epsilon_1) \supseteq \mathcal{D}_n(\theta, \epsilon_2) \supseteq \cdots
\end{equation}
Also, the feasible set $\Phi(u,\theta)$ is convex for fixed $u,\theta$ by \textbf{A1} and has a nonempty interior by \textbf{R1}.  This means $\Phi(u,\theta)$ is continuous in $\theta$ by Example 5.10 from \citep{rockafellar1998}, and so we can apply the Berge Maximum Theorem \citep{berge1963} to $\mathsf{FOP}$.  This implies $\mathcal{S}(u,\theta)$ is upper hemicontinuous in $\theta$ for fixed $u\in\mathcal{U}$.  By Remark 3.2 of \citep{dempe2015}, this means $Q_n(\theta)$ is lower semicontinuous.  Thus, by Proposition 7.4.d of \citep{rockafellar1998} we have that the extended real-valued function $\{Q_n(\theta; \epsilon_\nu)\ |\ \theta\in\Theta\}$ epiconverges to the extended real-valued function $\{Q_n(\theta)\ |\ \theta\in\Theta\}$.  The result then follows from Exercise 7.32.d and Theorem 7.33 of \citep{rockafellar1998}.
\Halmos\endproof

\begin{corollary}
\label{corollary:darqt}
Suppose $\mathbf{A1}$,$\mathbf{A2}$ and $\mathbf{R1}$ hold.  Given any $d > 0$, there exists $E,Z > 0$ such that if $\hat{\theta}_n \in z\text{-}\arg\min \{Q_n(\theta; \epsilon)\ |\ \theta\in\Theta\}$ for any $0 \leq z \leq Z$ and $0 \leq \epsilon \leq E$, then $\mathrm{dist}(\hat{\theta}_n,\arg\min \{Q_n(\theta)\ |\ \theta\in\Theta\}) < d$.
\end{corollary}

\proof{Proof. }
This is a restatement of Proposition \ref{prop:convmin}.
\Halmos\endproof

These results say that approximately solving $\mathsf{R}$--$\mathsf{IOP}$--$\mathsf{SAA}$ is equivalent to approximately solving $\mathsf{IOP}$--$\mathsf{SAA}$.

%

\subsection{Statistical Consistency}
\label{section:sc}

In order to prove statistical consistency, we will need to impose an additional regularity condition that ensures expectations of corresponding random variables exist.\\

\noindent\textbf{R2. } The set $\Theta$ is closed and bounded, and $\mathbb{E}(y^2) < +\infty$.\\

\noindent This regularity assumption ensures that the law of large numbers \citep{wald1949,jennrich1969,van2000} holds in our setting. The above expectation condition holds in many situations, including when $\mathcal{Y}$ is bounded or when $y$ has a sub-exponential distribution \citep{vershynin2012}.  This allows for settings where \textbf{IC} holds with measurement noise that is Gaussian, Bernoulli, bounded support, Laplacian, Exponential, and many other distributions. 


Our first statistical consistency result is that solving \textsf{R}--\textsf{IOP}--\textsf{SAA} is risk consistent.  To state the result, we must formally define the regularized version of the inverse optimization problem.  The regularized risk is
\leqnomode\begin{equation}
\tag*{$\mathsf{R}$--$\mathsf{RISK}$}
Q(\theta; \epsilon) = \mathbb{E}\Big(\min_{x\in\mathcal{S}(u,\theta; \epsilon)}\|y - x\|^2\Big),
\end{equation}\reqnomode
where $\mathcal{S}(u,\theta; \epsilon) = \{x\in\mathbb{R}^d : f(x,u,\theta) \leq V(u,\theta) + \epsilon,\ g(x,u,\theta) \leq \epsilon\}$ is the set of $\epsilon$-optimal solutions to \textsf{FOP}.  For given $\epsilon > 0$, we define the regularized inverse optimization problem to be
\leqnomode\begin{equation}
\tag*{$\mathsf{R}$--$\mathsf{IOP}$}
\min\{Q(\theta; \epsilon)\ |\ \theta\in\Theta\}.
\end{equation}\reqnomode
The first statistical consistency result specifically concerns nearly-optimal solutions of \textsf{R}--\textsf{IOP}--\textsf{SAA}.  We say that a sequence of solutions $\hat{\theta}_n$ is nearly-optimal for \textsf{R}--\textsf{IOP}--\textsf{SAA} with fixed $\epsilon > 0$ in probability if for any $\delta > 0$ we have
\begin{equation}
\lim_{n\rightarrow\infty}\mathbb{P}\Big(\text{dist}\big(\hat{\theta}_n,\arg\min\{Q_n(\theta; \epsilon)\ |\ \theta\in\Theta\}\big) > \delta\Big) = 0.
\end{equation}

\begin{theorem}
\label{theorem:wrc}
Suppose $\mathbf{A1}$,$\mathbf{A2}$ and $\mathbf{R1}$,$\mathbf{R2}$ hold.  Given any fixed $\epsilon > 0$, if $\hat{\theta}_n$ is nearly-optimal for $\mathsf{R}$--$\mathsf{IOP}$--$\mathsf{SAA}$ in probability, then we have $Q(\hat{\theta}_n; \epsilon) \stackrel{p}{\longrightarrow} \min \big\{Q(\theta; \epsilon)\ \big|\ \theta\in\Theta\big\}$. 
\end{theorem}

\proof{Proof. }Proposition \ref{prop:cont} gives continuity of $Q_n(\theta; \epsilon)$.  Thus, we can apply the uniform law of large numbers \citep{jennrich1969}, which gives
\begin{equation}
\label{eqn:ulln}
\textstyle\sup_{\theta\in\Theta}\big|Q_n(\theta; \epsilon)-Q(\theta; \epsilon)\big| \stackrel{p}{\longrightarrow} 0.
\end{equation}
Consider any $\theta_0\in\arg\min\{Q(\theta; \epsilon)\ |\ \theta\in\Theta\}$ and any $\theta_1\in\arg\min\{Q_n(\theta; \epsilon)\ |\ \theta\in\Theta\}$.  By assumption $Q_n(\theta_1; \epsilon) \leq Q_n(\theta_0; \epsilon)$, and so we have
\begin{equation}
Q(\hat{\theta}_n; \epsilon) + Q_n(\hat{\theta}_n; \epsilon) - Q(\hat{\theta}_n; \epsilon) + Q_n(\theta_1; \epsilon) - Q_n(\hat{\theta}_n; \epsilon) \leq Q(\theta_0; \epsilon) + Q_n(\theta_0; \epsilon) - Q(\theta_0; \epsilon).
\end{equation}
Rearranging terms gives
\begin{equation}
\label{eqn:rart}
Q(\hat{\theta}_n; \epsilon) - Q(\theta_0) \leq |Q_n(\hat{\theta}_n; \epsilon) - Q(\hat{\theta}_n; \epsilon)| + |Q_n(\theta_1; \epsilon) - Q_n(\hat{\theta}_n; \epsilon)| + |Q_n(\theta_0; \epsilon) - Q(\theta_0; \epsilon)|.
\end{equation}
Recall (i) $Q(\theta_0; \epsilon) \leq Q(\hat{\theta}_n; \epsilon)$ by definition of $\theta_0$, (ii) $Q_n(\theta; \epsilon)$ is continuous, and (iii) $\hat{\theta}_n$ is nearly-optimal for \textsf{R}--\textsf{IOP}--\textsf{SAA} in probability.  Thus, combining these facts with (\ref{eqn:ulln}) and (\ref{eqn:rart}) gives that $Q(\hat{\theta}_n; \epsilon) - Q(\theta_0; \epsilon)\stackrel{p}{\longrightarrow}0$.  This is the desired result.
\Halmos\endproof

This result says that if we choose any $\epsilon > 0$ and solve $\mathsf{R}$--$\mathsf{IOP}$--$\mathsf{SAA}$ to generate an estimate $\hat{\theta}_n$, then the predictions given by the $\epsilon$-optimal solutions to \textsf{FOP} (i.e., $\mathcal{S}(u,\hat{\theta}_n;\epsilon)$) are asymptotically the best possible set of predictions when the error of predictions is measured using \textsf{R}-\textsf{RISK}.  A stronger risk consistency result is not possible in the general setting because $Q(\theta)$ is typically discontinuous, and so the above result can be interpreted as a weak consistency result.  

A stronger risk consistency result is possible in the case where $f(x,u,\theta)$ is strictly convex.  We say that a sequence of solutions $\hat{\theta}_n$ is nearly-optimal for \textsf{IOP}--\textsf{SAA} in probability if for any $\delta > 0$ we have\footnote{}\footnotetext[2]{Note that this notion of near-optimality is defined with respect to $\textsf{IOP}$--$\textsf{SAA}$, whereas the definition of near-optimality given in (15) is with respect to the regularized formulation $\textsf{R}$--$\textsf{IOP}$--$\textsf{SAA}$.} 
\begin{equation}
\lim_{n\rightarrow\infty}\mathbb{P}\Big(\text{dist}\big(\hat{\theta}_n,\arg\min\{Q_n(\theta)\ |\ \theta\in\Theta\}\big) > \delta\Big) = 0.
\end{equation} 
\begin{theorem}
Suppose $\mathbf{A1}$,$\mathbf{A2}$ and $\mathbf{R1}$,$\mathbf{R2}$ hold.  If $f(x,u,\theta)$ is strictly convex in $x$ (for fixed $u\in\mathcal{U}$ and $\theta\in\Theta$) and $\hat{\theta}_n$ is nearly-optimal for $\mathsf{IOP}$--$\mathsf{SAA}$ in probability, then we have $Q(\hat{\theta}_n) \stackrel{p}{\longrightarrow} \min \big\{Q(\theta)\ \big|\ \theta\in\Theta\big\}$.
\end{theorem}

\proof{Proof. }Proposition \ref{proposition:sccont} gives continuity of $Q_n(\theta)$.  The remainder of the proof is identical to Theorem \ref{theorem:wrc}.
\Halmos\endproof

This result says that when \textsf{FOP} is a strictly convex optimization problem and we solve $\mathsf{IOP}$--$\mathsf{SAA}$ to generate an estimate $\hat{\theta}_n$, then the predictions given by the solutions to \textsf{FOP} (i.e., $\mathcal{S}(u,\hat{\theta}_n)$) are asymptotically the best possible set of predictions when the error of predictions is measured using \textsf{RISK}.  The reason it is possible to show risk consistency in this case is that $Q(\theta)$ will be continuous in this setting.

Our final statistical consistency result is that solving \textsf{IOP}--\textsf{SAA} is estimation consistent when \textbf{IC} holds.

\begin{theorem}
\label{theorem:estconst}
Suppose $\mathbf{A1}$,$\mathbf{A2}$ and $\mathbf{R1}$,$\mathbf{R2}$ and $\mathbf{IC}$ hold.  If $\hat{\theta}_n$ is nearly-optimal for $\mathsf{IOP}$--$\mathsf{SAA}$ in probability, then we have $\hat{\theta}_n \stackrel{p}{\longrightarrow} \theta_0$. 
\end{theorem}

\proof{Proof. }Because the feasible set $\Phi(u,\theta)$ is convex for fixed $u,\theta$ by \textbf{A1} and has a nonempty interior by \textbf{R1}, this means $\Phi(u,\theta)$ is continuous in $\theta$ by Example 5.10 from \citep{rockafellar1998}.  Thus, we can apply the Berge Maximum Theorem \citep{berge1963} to $\mathsf{FOP}$.  This implies $\mathcal{S}(u,\theta)$ is upper hemicontinuous in $\theta$ for fixed $u\in\mathcal{U}$.  By Remark 3.2 of \cite{dempe2015}, this means $Q_n(\theta)$ is lower semicontinuous.  Thus, we can apply Theorem 5.14 of \citep{van2000}.\footnote{}\footnotetext[3]{Technically, this theorem applies to maximizing upper semicontinuous functions, but the results and proof trivially extend to the case of minimizing lower semicontinuous functions.}  The result follows from the conclusion of Theorem 5.14 of \citep{van2000} if we can show (i) $\theta_0 \in \arg\min\{Q(\theta)\ |\ \theta\in\Theta\}$, and that (ii) $\theta_0$ is the unique solution.  First, note $Q(\theta) = \mathbb{E}(\min_{x\in\mathcal{S}(u,\theta)}\|\xi - x\|^2) + \mathbb{E}(w^2)$, since $\xi,x$ is almost surely independent of $w$ because by \textbf{IC} we have that (i) $\xi,u$ are independent of $w$, and (ii) $\mathcal{S}(u,\theta)$ is almost surely single-valued.  Since by \textbf{IC} we have $\xi \in \mathcal{S}(u,\theta_0)$, this means that $Q(\theta_0) = \mathbb{E}(w^2)$ and that $\theta_0 \in \arg\min\{Q(\theta)\ |\ \theta\in\Theta\}$.  Next, consider any $\theta\in\Theta\setminus\theta_0$.  Then by \textbf{IC} we have $\mathbb{E}[\min_{x\in\mathcal{S}(u,\theta)}\|\xi - x\|^2\ |\ u\in\mathcal{U}(\theta)] > 0$ since $\xi \in \mathcal{S}(u,\theta_0)$ and $\text{dist}(\mathcal{S}(u,\theta), \mathcal{S}(u,\theta_0)) > 0$ for each $u\in\mathcal{U}(\theta)$.  Because $\mathbb{P}(u \in \mathcal{U}(\theta)) > 0$ from \textbf{IC}, this means $\mathbb{E}(\min_{x\in\mathcal{S}(u,\theta)}\|\xi - x\|^2) > 0$ for any $\theta\in\Theta\setminus\theta_0$.  Consequently, we have $Q(\theta) > Q(\theta_0)$ for any $\theta\in\Theta\setminus\theta_0$.
\Halmos\endproof

\section{Numerical Approaches to Solving \textsf{IOP}--\textsf{SAA}}
\label{section:numerical}
Solving $\mathsf{IOP}$--$\mathsf{SAA}$ with $Q_n(\theta)$ as formulated in $\mathsf{DB}$--$\mathsf{RISK}$--$\mathsf{SAA}$ is still difficult because it is a nonconvex problem even under \textbf{A1},\textbf{A2},\textbf{R1}.  We will propose two approaches to solving this problem.  The first is an enumeration algorithm that is applicable to situations where $p$ is modest (i.e., the $\theta \in \mathbb{R}^p$ parameter has between 1 to 5 dimensions).  The second approach we describe is a semiparametric algorithm, and it can be used in cases where $\theta \in \mathbb{R}^p$ is higher-dimensional and the noise term $w$ has a specific distribution. For both algorithms, we will prove that the estimates computed by these methods satisfy the conditions required for statistical consistency.

The difference in the two algorithms is how they trade-off computational and statistical performance.  The enumeration algorithm requires computation exponential in $p$, while the semiparametric algorithm needs computation polynomial $p$ computation. But the statistical performance of the methods will be the opposite. The estimates and risk of the enumeration algorithm are anticipated to converge at faster rate (with respect to the number of data points) than those of the semiparametric algorithm. The reason is that the semiparametric algorithm makes use of a nonparametric step (via the L2NW estimator), which is well-known to generally converge at a slower rate than a fully parametric approach.  Precisely characterizing the statistical convergence rates of the two algorithms is left open for future work.


Though the enumeration algorithm needs exponential in $p$ computation, it is still practical for many real-world problems. Many principal-agent problems (e.g. \cite{zhang2008dynamic,crama2008milestone}) use models where the parameter set is modest in dimensionality (i.e., utility functions with 2 or 3 \emph{type} parameters).  We demonstrate the practicality of the enumeration algorithm in Section 5 through an energy-related example using real data.

\subsection{Enumeration Algorithm}

The main idea of this algorithm is that computing $Q_n(\theta)$ and $Q_n(\theta; \epsilon)$ for fixed values of $\theta$ can be done in polynomial time since $\mathsf{DB}$--$\mathsf{RISK}$--$\mathsf{SAA}$ and $\mathsf{R}$--$\mathsf{DB}$--$\mathsf{RISK}$--$\mathsf{SAA}$ are convex optimization problems by Propositions \ref{prop:convex} and \ref{prop:convex_r}, respectively.  This approach enumerates over different fixed values of $\theta$ and solves a series of polynomial time problems.  However, $\Theta$ is a continuous set since because it is convex by \textbf{A2}.  To enable enumeration, we discretize $\Theta$ using a $\delta$-net of $\Theta$, which we will call $\mathcal{T}(\delta)$.  (Here, we define this to mean that $\mathcal{T}(\delta)$ is a finite set such that $\max_{\theta\in\Theta}\min_{t \in \mathcal{T}(\delta)} \|t - \theta\| \leq \delta$.)  We then compute $Q_n(\theta; \epsilon)$ for all $\theta\in\mathcal{T}(\delta)$.  And our approximate solution is finally given by $\hat{\theta}_n = \arg\min\{Q_n(\theta; \epsilon)\ |\ \theta \in\mathcal{T}(\delta)\}$.  

This approach requires continuity of $Q_n(\theta; \epsilon)$ because otherwise performing an enumeration via the $\delta$-net $\mathcal{T}(\delta)$ may not get sufficiently close to the optimal value.  However, $Q_n(\theta; \epsilon)$ is only guaranteed to be continuous at $\epsilon = 0$ when $f(x,u,\theta)$ is strictly convex for fixed $u,\theta$ by Proposition \ref{proposition:sccont} and since $Q_n(\theta; 0) = Q_n(\theta)$ by definition.  Hence, we require $\epsilon > 0$ for cases where $f(x,u,\theta)$ is \emph{not} strictly convex to ensure continuity of $Q_n(\theta;\epsilon)$ by Proposition \ref{prop:cont}.  Of course, when $f(x,u,\theta)$ is strictly convex we can set $\epsilon = 0$ and maintain continuity of $Q_n(\theta; \epsilon)$.

This approach is formally presented in Algorithm \ref{algorithm:enum}. Importantly, it can be shown that this enumeration algorithm generates nearly-optimal solutions of \textsf{IOP}--\textsf{SAA} and \textsf{R}--\textsf{IOP}--\textsf{SAA}.  This means the solutions computed by this algorithm satisfy the conditions in Theorems 2, 3, and 4 that are needed for statistical consistency.  In practice, $\epsilon$ is chosen to be $\epsilon = 0$ when \textsf{FOP} is strictly convex, and otherwise $\epsilon$ is chosen to be a small positive value that controls the desired precision of the resulting estimate. An appropriate approach to choose $\epsilon$ and $\delta$ is to use cross-validation, which is a standard data-driven approach from statistics for choosing such parameters \citep{hastie2009}.

\begin{figure}
\begin{algorithm}[H]
 \KwData{fixed $\delta> 0$ and $\epsilon \geq 0$}
 \KwResult{estimate $\hat{\theta}_n$}
 set $\mathcal{T}(\delta)$ to be $\delta$-net of $\Theta$\;
\ForEach{$\theta \in \mathcal{T}(\delta)$}{
compute $Q_n(\theta; \epsilon)$ by solving $\mathsf{R}$--$\mathsf{DB}$--$\mathsf{RISK}$--$\mathsf{SAA}$\;
}
set $\hat{\theta}_n \in \arg\min\{Q_n(\theta; \epsilon)\ |\ \theta \in\mathcal{T}(\delta)\}$\;
 \caption{Enumeration Algorithm}
\label{algorithm:enum}
\end{algorithm}
\end{figure}

\begin{theorem}
Suppose $\mathbf{A1}$,$\mathbf{A2}$ and $\mathbf{R1}$ hold.  Given any $d > 0$, there exists $E,\Delta > 0$ such that if $\hat{\theta}_n$ is computed using the enumeration algorithm (i.e., Algorithm \ref{algorithm:enum}) for any $0 < \epsilon \leq E$ and $0 < \delta \leq \Delta$, then $\mathrm{dist}(\hat{\theta}_n,\arg\min \{Q_n(\theta)\ |\ \theta\in\Theta\}) < d$.
\end{theorem}

\proof{Proof. }By Corollary \ref{corollary:darqt}, there exists $E,Z > 0$ such that if $\hat{\theta}_n \in z\text{-}\arg\min \{Q_n(\theta; \epsilon)\ |\ \theta\in\Theta\}$ for any $0 \leq z \leq Z$ and $0 \leq \epsilon \leq E$, then $\mathrm{dist}(\hat{\theta}_n,\arg\min \{Q_n(\theta)\ |\ \theta\in\Theta\}) < d$.  Suppose we choose $z = Z$.  Because $Q_{n}(\theta; \epsilon)$ is continuous in $\theta$ by Proposition \ref{prop:cont}, there exists $\Delta > 0$ such that for any $0 < \delta \leq \Delta$ we have
\begin{equation}
\min\big\{Q_{n}(\theta; \epsilon) -  Q_n(\theta_0; \epsilon)\ \big|\ \theta \in \mathcal{T}(\delta)\big\} < z,
\end{equation}
where $\theta_0 \in \arg\min\{Q_n(\theta; \epsilon)\ |\ \theta \in\Theta\}$.  By construction, we have
\begin{equation}
\arg\min\{Q_n(\theta; \epsilon)\ |\ \theta\in\mathcal{T}(\delta)\} \subseteq z\text{-}\arg\min \{Q_n(\theta; \epsilon)\ |\ \theta\in\Theta\}.
\end{equation}
Next, note the enumeration algorithm returns a solution $\hat{\theta}_n \in \arg\min\{Q_n(\theta; \epsilon)\ |\ \theta\in\mathcal{T}(\delta)\}$, which also satisfies $\hat{\theta}_n \in z\text{-}\arg\min \{Q_n(\theta; \epsilon)\ |\ \theta\in\Theta\}$.  The result follows from applying the first line of the proof.
\Halmos\endproof
Theorem 5 states that the estimate obtained using the enumeration algorithm will be at most a distance of $d$ from the set of optimal solutions to \textsf{IOP-SAA}. It immediately follows that for small $d$, the solution of the enumeration algorithm will retain the desirable statistical properties of the solutions to \textsf{IOP-SAA}. As mentioned above, in the special case where \textsf{FOP} is a strictly convex optimization problem we can simplify the algorithm by setting $\epsilon = 0$.  We have a corresponding result about the correctness of the algorithm in this case.

\begin{theorem}
Suppose $\mathbf{A1}$,$\mathbf{A2}$ and $\mathbf{R1}$ hold.  If $f(x,u,\theta)$ is strictly convex in $x$ (for fixed $u\in\mathcal{U}$ and $\theta\in\Theta$), then given any $d > 0$ there exists $\Delta > 0$ such that if $\hat{\theta}_n$ is computed using the enumeration algorithm for $\epsilon = 0$ and any $0 < \delta \leq \Delta$, then $\mathrm{dist}(\hat{\theta}_n,\arg\min \{Q_n(\theta)\ |\ \theta\in\Theta\}) < d$.
\end{theorem}

\proof{Proof. }By Corollary \ref{corollary:darqt}, there exists $E,Z > 0$ such that if $\hat{\theta}_n \in z\text{-}\arg\min \{Q_n(\theta; \epsilon)\ |\ \theta\in\Theta\}$ for any $0 \leq z \leq Z$ and $0 \leq \epsilon \leq E$, then $\mathrm{dist}(\hat{\theta}_n,\arg\min \{Q_n(\theta)\ |\ \theta\in\Theta\}) < d$.  Suppose we choose $z = Z$ and $\epsilon = 0$, and note that $Q_n(\theta; 0) = Q_n(\theta)$ by their definitions.  Because $Q_{n}(\theta)$ is continuous in $\theta$ by Proposition \ref{proposition:sccont}, there exists $\Delta > 0$ such that for any $0 < \delta \leq \Delta$ we have
\begin{equation}
\min\big\{Q_{n}(\theta) -  Q_n(\theta_0)\ \big|\ \theta \in \mathcal{T}(\delta)\big\} < z,
\end{equation}
where $\theta_0 \in \arg\min\{Q_n(\theta)\ |\ \theta \in\Theta\}$.  By construction, we have
\begin{equation}
\arg\min\{Q_n(\theta)\ |\ \theta\in\mathcal{T}(\delta)\} \subseteq z\text{-}\arg\min \{Q_n(\theta)\ |\ \theta\in\Theta\}.
\end{equation}
Next, note the enumeration algorithm returns a solution $\hat{\theta}_n \in \arg\min\{Q_n(\theta)\ |\ \theta\in\mathcal{T}(\delta)\}$, which also satisfies $\hat{\theta}_n \in z\text{-}\arg\min \{Q_n(\theta)\ |\ \theta\in\Theta\}$.  The result follows from the first line of the proof.
\Halmos\endproof

\subsection{Semiparametric Approach}

Our second approach to solving $\mathsf{IOP}$--$\mathsf{SAA}$ is a semiparametric approach. We will need to make an additional assumption about the structure of the problem, as well as impose two more regularity conditions, in order to be able use this approach.  We begin with the additional assumption.\\

\noindent\textbf{A3. } The constraint function $g(x,u,\theta)$ is independent of $\theta$, meaning it can be written as $g(x,u,\theta) = g_0(x,u)$.  The objective function $f(x,u,\theta)$ is affine in $\theta$, meaning it can be written as
\begin{equation}
f(x,u,\theta) = f_0(x,u) + \sum_{j=1}^p  \theta_j f_j(x,u).
\end{equation}
Independence of the constraint $g$ from $\theta$ is required because the semiparametric approach relies on fully knowing the feasible region of the forward problem. We note that this is not a particularly strong assumption, since in utility estimation settings one would expect the unknown parameters to appear in the objective function of the forward problem. \cite{keshavarz2011} and \cite{bertsimas2013} also assume that the feasible region of the forward problem is independent of the unknown parameters.  The second part of ${\bf A3}$ ensures that the Lagrangian dual function $h(\lambda,u,\theta)$ is concave in $\theta$. This will enable efficient computation in our semiparametric approach.  Next, we describe the two additional regularity conditions.  The first is\\

\noindent\textbf{R3. } The objective function $f(x,u,\theta)$ is strictly convex in $x$ (for fixed $u \in\mathcal{U}$ and $\theta \in\Theta$) and twice continuously differentiable in $x,u,\theta$, and the constraints $g(x,u,\theta)$ are continuously differentiable in $x,u,\theta$.\\

\noindent Condition \textbf{R3} ensures smoothness in the objective function and constraints. The reason we also include a strict convexity assumption is that it acts as a regularity condition: Strictly speaking, we require uniqueness of solutions to \textsf{FOP} (which is needed for the de-noising step in our semiparametric algorithm) and a second-order growth condition
\begin{equation}
f(x,u,\theta) \geq V(u,\theta) + c\cdot [\text{dist}(x,\mathcal{S}(u,\theta))]^2,
\end{equation}
for some $c > 0$ and all $x \in \Phi(u,\theta)$ (which ensures H\"{o}lder continuity of the solution set $\mathcal{S}(u,\theta)$ with degree $1/2$ \citep{bonnans2000}). Unfortunately, this growth condition can be difficult to directly check even though it has been completely characterized for convex optimization problems \citep{bonnans1995}.  Fortunately, strict convexity with Slater's constraint qualification (which holds under {\bf R1}) implies both uniqueness of solutions to \textsf{FOP} and this second-order growth condition \citep{bonnans1995second}.  Hence {\bf R3} is sufficient for proving statistical convergence using our algorithm. We also note that our results could be extended to the case where the problem satisfies the first-order growth condition
\begin{equation}
f(x,u,\theta) \geq V(u,\theta) + c\cdot \text{dist}(x,\mathcal{S}(u,\theta)),
\end{equation}
for some $c > 0$ and all $x \in \Phi(u,\theta)$. Under this alternate growth condition, the solution set is H\"{o}lder continuous with degree $1$ (instead of $1/2$). This affects the bound expression in Proposition \ref{prop:unifconv} slightly, but otherwise does not qualitatively change our results.
\\

\noindent\textbf{R4. } The noise random variable $w$ has a sub-exponential distribution, meaning there exists $c > 0$ such that $\mathbb{P}(|w| > t) \leq \exp(1-t/c)$.  Also, the probability density function $\mu(u)$ of $u$ is continuously differentiable and is bounded from zero (i.e., $\min_{u\in\mathcal{U}}\mu(u) > 0$).\\

\noindent This regularity condition ensures the distribution of the random variables $w,u$ are not extreme.  Most commonly used heavy-tailed noise distributions are sub-exponential distributions, and so $\textbf{R4}$ is satisfied by Gaussian, Bernoulli, bounded support, Laplacian, Exponential, and many other distributions \citep{vershynin2012}.  Also, the regularity condition on $\mu(u)$ implies $\mathcal{U}$ is bounded.

The idea behind the semiparametric approach is the observation that $\mathsf{R}$--$\mathsf{DB}$--$\mathsf{RISK}$--$\mathsf{SAA}$ is convex in $\theta$ for fixed $x$ when \textbf{A3} holds.  However, because the $y_i$ are measured with noise, we cannot simply make the substitution $x_i = y_i$.  To overcome this difficulty, we first de-noise the $y_i$ using a nonparametric estimator.  Specifically, we define the $\ell_2$-regularized Nadaraya-Watson (L2NW) estimator \citep{aswani2013_automatica} as
\begin{equation}
\label{eqn:nw}
\overline{x}_i = \frac{\gamma^{-m}\cdot\frac{1}{n}\sum_{j=1}^n y_j\cdot K\big(\frac{u_j-u_i}{\gamma}\big)}{\sigma + \gamma^{-m}\cdot\frac{1}{n}\sum_{j=1}^nK\big(\frac{u_j-u_i}{\gamma}\big)}, 
\end{equation}
where $\gamma > 0$ is the \emph{bandwidth} parameter, $\sigma > 0$ is the $\ell_2$-\emph{regularization} parameter, and $K : \mathbb{R}^m\rightarrow\mathbb{R}$ is a \emph{kernel function} that satisfies the following properties (i) $K(u) \geq 0$, (ii) $K(u) = 0$ for $\|u\| > 1$, (iii) $K(u) = K(-u)$, and (iv) $\int K(u)du = 1$.  A common example of a kernel function is the Epanechnikov kernel, which is defined as the function
\begin{equation}
K(u) = \begin{cases} \textstyle\frac{3}{4}\cdot(1-\|u\|^2), & \text{if } \|u\| \leq 1\\
0, &\text{otherwise}\end{cases}
\end{equation}
The L2NW estimator (\ref{eqn:nw}) is computed in polynomial time, and it serves to de-noise the $x_i$ in the manner described by the following proposition.


\begin{proposition}
\label{prop:unifconv}
Suppose $\mathbf{A1}$ and $\mathbf{R1}$--$\mathbf{R4}$ hold.  If $\gamma = O(n^{-2/(8m+1)})$ and $\sigma = O(\gamma)$, then $\mathcal{S}(u,\theta)$ consists of a single point, and for sufficiently large $n$ we have
we have 
\begin{equation}
\mathbb{P}\Big(\max_{i\in[n]}\big\|\overline{x}_i - \mathcal{S}(u_i,\theta_0)\big\| > n^{-1/(18m)}\Big) \leq k_1\exp\Big(-k_2n^{1/4}\Big),
\end{equation}
where $k_1,k_2 > 0$ are constants.  In particular, this implies $\max_{i\in[n]}\big\|\overline{x}_i - \mathcal{S}(u_i,\theta_0)\big\|\stackrel{p}{\longrightarrow} 0$.
\end{proposition}

\proof{Proof. }The first part follows from the strict convexity assumption in \textbf{R3}, and the third part follows directly from the second part.  And so we focus on proving the second part.  We will prove this using a truncation argument (see for instance \citep{tao2012}).


First, note that the function $\psi(x,y) = x/y$ over the domain $(x,y) \in [-M,M]\times[\sigma,\sigma+1]$ is Lipschitz continuous with constant $L_1 = \sqrt{(M^2+(\sigma+1)^2)}/\sigma^2$.  Suppose we choose $M = \max_{u\in\mathcal{U}}\|\mu(u)\mathcal{S}(u,\theta_0)\| + 1$.  As a result, using Lemma \ref{lemma:1} and Lemma \ref{lemma:2} we have
\begin{equation}
\label{eqn:1}
\begin{aligned}
&\textstyle\mathbb{P}\Big(\big\|\overline{x}_i - \frac{\mu(u_i)\mathcal{S}(u_i,\theta_0)}{\sigma + \mu(u_i)}\big\| > t\Big) \\
&\quad\leq \textstyle\mathbb{P}\Big(\big|\gamma^{-m}\cdot\frac{1}{n}\sum_{j=1}^nK\big(\frac{u_j-u_i}{\gamma}\big) - \mu(u_i)\big|>t/L_1\Big) + \textstyle\mathbb{P}\Big(\big\|\gamma^{-m}\cdot\frac{1}{n}\sum_{j=1}^ny_j\cdot K\big(\frac{u_j-u_i}{\gamma}\big)\big\|>M\Big)+\\
&\qquad\qquad\textstyle\mathbb{P}\Big(\big\|\gamma^{-m}\cdot\frac{1}{n}\sum_{j=1}^ny_j\cdot K\big(\frac{u_j-u_i}{\gamma}\big) - \mu(u_i)\mathcal{S}(u_i,\theta_0)\big\|>t/L_1\Big)\\
&\quad\leq 2\exp\Big(-2c_2 n\gamma^{2m}\cdot(t/L_1-c_1\cdot \gamma)^2\Big) + 2\exp\Big(-2c_2 n\gamma^{2m}\cdot(1-c_1\cdot \gamma)^2\Big)+\\
&\qquad\qquad2\exp\Big(-2c_5 n\gamma^{2m}\cdot(t/L_1-c_3\cdot \gamma^{1/2}-c_4\cdot \gamma)\Big),
\end{aligned}
\end{equation}
for $t > \max\{c_1\cdot \gamma,c_3\cdot \gamma^{1/2}+c_4\cdot \gamma\}$.
Next, observe that the function $\psi(x,y)$ over the domain
\begin{equation}
(x,y) \in [\min_{u\in\mathcal{U}}\mu(u)\mathcal{S}(u,\theta),\max_{u\in\mathcal{U}}\mu(u)\mathcal{S}(u,\theta)]\times[\min_{u\in\mathcal{U}}\mu(u),\max_{u\in\mathcal{U}}\mu(u)],
\end{equation}
is Lipschitz continuous with some constant $L_2 > 0$ since (i) the denominator of $\psi$ is bounded away from zero because of \textbf{R4}, and (ii) the numerator of $\psi$ is bounded by \textbf{R1},\textbf{R4}.  Thus, we have
\begin{equation}
\label{eqn:2}
\mathbb{P}\Big(\big\|\overline{x}_i - \mathcal{S}(u_i,\theta_0)\big\| > t\Big) \leq \mathbb{P}\Big(\big\|\overline{x}_i - \textstyle\frac{\mu(u_i)\mathcal{S}(u_i,\theta_0)}{\sigma + \mu(u_i)}\big\| > t-\sigma/L_2\Big),
\end{equation}
for $t > \sigma/L_2$.  Suppose we choose $\gamma = O(n^{-2/(8m+1)})$, $\sigma = O(\gamma)$, and $t = n^{-1/(16m+2)}$.  Then combining (\ref{eqn:1}) and (\ref{eqn:2}) gives that for sufficiently large $n$ we have
\begin{equation}
\label{eqn:3}
\mathbb{P}\Big(\big\|\overline{x}_i - \mathcal{S}(u_i,\theta_0)\big\| > n^{-1/(16m+2)}\Big) \leq c_6\exp\Big(-c_7n^{1/2}\Big),
\end{equation}
where $c_6,c_7 > 0$ are constants.  And so combining the union bound with (\ref{eqn:3}) gives
\begin{equation}
\begin{aligned}
\mathbb{P}\Big(\max_{i\in[n]}\big\|\overline{x}_i - \mathcal{S}(u_i,\theta_0)\big\| > n^{-1/(16m+2)}\Big) &\leq n\mathbb{P}\Big(\big\|\overline{x}_i - \mathcal{S}(u_i,\theta_0)\big\| > n^{-1/(16m+2)}\Big)\\
&\leq c_6\exp\Big(-c_7n^{1/2} + \log n\Big).
\end{aligned}
\end{equation}
The final implication of the result follows by noting that $n^{-2/(8m+1)}\rightarrow  0$ and $c_6\exp(-c_7n^{1/2} + \log n) \rightarrow 0$ as $n\rightarrow\infty$.
\Halmos\endproof

Before we present our algorithm, we need one more result that provides additional understanding for the semiparametric approach.  Consider the following optimization problem
\leqnomode\begin{equation}
\tag*{$\mathsf{ROBUST}$--$\mathsf{IOP}$--$\mathsf{SAA}$}
\min_\theta\max_{\epsilon\geq0}\big\{Q_n(\theta; \epsilon)\ \big|\ \theta\in\Theta\big\},
\end{equation}\reqnomode

\begin{proposition}
Suppose $\mathbf{A1}$,$\mathbf{A2}$ and $\mathbf{R1}$ hold.  Then the solution sets in $\theta$ of $\mathsf{ROBUST}$--$\mathsf{IOP}$--$\mathsf{SAA}$ and $\mathsf{IOP}$--$\mathsf{SAA}$ are equivalent, and the optimal value of $\mathsf{ROBUST}$--$\mathsf{IOP}$--$\mathsf{SAA}$ occurs at $\epsilon = 0$.
\end{proposition}

\proof{Proof. }Let $\mathcal{C}_n(\theta, \epsilon)$ be the feasible set of $\mathsf{R}$--$\mathsf{DB}$--$\mathsf{RISK}$--$\mathsf{SAA}$.  As shown in the proof for Proposition \ref{prop:convmin}, the feasible set satisfies
\begin{equation}
\mathcal{C}_n(\theta, 0) \subseteq \mathcal{C}_n(\theta, \epsilon),
\end{equation}
for all $\epsilon \geq 0$.  As a result, we must have that $Q_n(\theta; 0) \geq Q_n(\theta; \epsilon)$ for all $\epsilon \geq 0$.  This means that $\max_{\epsilon\geq0}Q_n(\theta; \epsilon) = Q_n(\theta; 0)$.  The result holds because $Q_n(\theta; 0) = Q_n(\theta)$ by definition.
\Halmos\endproof

Given the above relationship that the optimal value of $\mathsf{ROBUST}$--$\mathsf{IOP}$--$\mathsf{SAA}$ occurs at $\epsilon = 0$, we propose to solve the inverse optimization problem using the following formulation:
\leqnomode\begin{equation}
\tag*{$\mathsf{SP}$--$\mathsf{IOP}$--$\mathsf{RISK}$--$\mathsf{SAA}$}
\begin{aligned}
\hat{\theta}_n \in \arg\min\ &\frac{1}{n}\sum_{i=1}^n\epsilon_i\\
\text{s.t. }&f(\overline{x}_i,u_i,\theta)-h(\lambda_i,u_i,\theta)\leq \epsilon_i,&\forall i\in[n]\\
&\lambda_i\geq0,&\forall i\in[n]
\end{aligned}
\end{equation}\reqnomode
where the $\overline{x}_i$ are as defined in (\ref{eqn:nw}).  This is a convex optimization problem.

\begin{proposition}
Suppose $\mathbf{A1}$--$\mathbf{A3}$ and $\mathbf{R1}$ hold.  Then $\mathsf{SP}$--$\mathsf{IOP}$--$\mathsf{RISK}$--$\mathsf{SAA}$ is a convex optimization problem.
\end{proposition}

We now have the elements to construct our semiparametric algorithm, which is a two-step approach.  In the first step, we de-noise the $y_i$ data using the L2NW estimator given in (\ref{eqn:nw}). This de-noising step produces an estimate of the true underlying optimal solution, which we represent by $\bar{x}_i$. While the estimates $\bar{x}_i$ are asymptotically (in $n$) optimal (cf. Proposition \ref{prop:unifconv}), they may be suboptimal at finite $n$. Therefore, in the second step, we solve $\mathsf{SP}$--$\mathsf{IOP}$--$\mathsf{RISK}$--$\mathsf{SAA}$, which produces a parameter estimate $\hat{\theta}_n$ that minimizes the suboptimality of $\bar{x}_i$. This approach maintains statistical consistency because the $\bar{x}_i$ are denoised, and it is formally presented in Algorithm \ref{algorithm:sp}.  Importantly, it can be shown that this semiparametric algorithm generates nearly-optimal solutions of \textsf{IOP}--\textsf{SAA}. This means the solutions computed by this algorithm satisfy the conditions in Theorems 2, 3, and 4 that are needed for statistical consistency.   In practice, the values of $\sigma$ and $\gamma$ can be can be chosen using cross-validation, which is a standard data-driven approach from statistics for choosing such parameters \citep{hastie2009}.

\begin{figure}
\begin{algorithm}[H]
 \KwData{fixed $\gamma > 0$ and $\sigma > 0$}
 \KwResult{estimate $\hat{\theta}_n$}
\ForEach{$i \in [n]$}{
compute $\overline{x}_i$ using using (\ref{eqn:nw})\;
}
compute $\hat{\theta}_n$ using $\mathsf{SP}$--$\mathsf{IOP}$--$\mathsf{RISK}$--$\mathsf{SAA}$\;
 \caption{Semiparametric Algorithm}
\label{algorithm:sp}
\end{algorithm}
\end{figure}
\begin{theorem}
Suppose $\mathbf{A1}$--$\mathbf{A3}$ and $\mathbf{R1}$--$\mathbf{R4}$ and $\mathbf{IC}$ hold.  If $\sigma = O(n^{-2/(8m+1)})$, $\lambda = O(\sigma)$, and $\hat{\theta}_n$ is computed using the semiparametric algorithm (i.e., Algorithm \ref{algorithm:sp}) ; then $\hat{\theta}_n$ is nearly-optimal for $\mathsf{IOP}$--$\mathsf{SAA}$ in probability.
\end{theorem}

\proof{Proof. }Note that $\min \{-h(\lambda,u,\theta)\ | \lambda \geq 0\} = -f(\mathcal{S}(u,\theta),u,\theta)$ by strong duality (which holds because of \textbf{A1},\textbf{R1} \citep{bonnans2000}).  Next, consider the function
\begin{equation}
R(\theta) = \mathbb{E}\Big(\min_{\lambda \geq 0}f(\mathcal{S}(u,\theta_0),u,\theta)-h(\lambda,u,\theta)\Big) = \mathbb{E}\Big(f(\mathcal{S}(u,\theta_0),u,\theta)-f(\mathcal{S}(u,\theta),u,\theta)\Big),
\end{equation}
its sample average approximation
\begin{equation}
R_n(\theta) = \frac{1}{n}\sum_{i=1}^n\Big(\min_{\lambda_i \geq 0}f(\mathcal{S}(u_i,\theta_0),u_i,\theta)-h(\lambda_i,u_i,\theta)\Big) = \frac{1}{n}\sum_{i=1}^n\Big(f(\mathcal{S}(u_i,\theta_0),u_i,\theta)-f(\mathcal{S}(u_i,\theta),u_i,\theta)\Big),
\end{equation}
and its semiparametric approximation
\begin{equation}
\label{eqn:seap}
\overline{R}_n(\theta) = \frac{1}{n}\sum_{i=1}^n\Big(\min_{\lambda_i \geq 0}f(\overline{x}_i,u_i,\theta)-h(\lambda_i,u_i,\theta)\Big) = \frac{1}{n}\sum_{i=1}^n\Big(f(\overline{x}_i,u_i,\theta)-f(\mathcal{S}(u_i,\theta),u_i,\theta)\Big).
\end{equation}
Note that $\min\{\overline{R}_n(\theta)\ |\ \theta\in\Theta\}$ is simply a reformulation of $\mathsf{SP}$--$\mathsf{IOP}$--$\mathsf{RISK}$--$\mathsf{SAA}$.  Next, observe that $\mathbb{E}[f(\mathcal{S}(u,\theta_0),u,\theta)-f(\mathcal{S}(u,\theta),u,\theta)|u\in\mathcal{U}(\theta)] > 0$ since (i) $f(x,u,\theta)$ is twice continuously differentiable in $x$ by \textbf{R3}, and (ii) $\text{dist}(\mathcal{S}(u,\theta), \mathcal{S}(u,\theta_0)) > 0$ for each $u\in\mathcal{U}(\theta)$ by \textbf{IC}.  Consequently, we have $R(\theta) > 0$ for $\theta\in\Theta\setminus\theta_0$.  As shown in the proof for Proposition \ref{proposition:sccont}, $\mathcal{S}(u,\theta)$ is continuous in $\theta$.  And so $R_n(\theta)$ and $\overline{R}_n(\theta)$ are continuous because (i) $f(x,u,\theta)$ is twice continuously differentiable in $x,\theta$ by \textbf{R3}.  

Next, recall that $\mathcal{U}$ is bounded by \textbf{R4}, $\Theta$ is bounded by \textbf{R2}, $f(x,u,\theta)$ is twice continuously differentiable in $x,\theta$ by \textbf{R3}, and the feasible set of \textsf{FOP} is absolutely bounded by \textbf{R1}.  This means there exists $L > 0$ such that for all $\theta\in\Theta$ we have $\max_{i\in[n]}|f(\overline{x},u_i,\theta) - f(\mathcal{S}(u_i,\theta_0),u_i,\theta)| \leq Ln^{-1/(18m)}$ whenever $\max_{i\in[n]}\|\overline{x}_i-\mathcal{S}(u_i,\theta_0)\| \leq n^{-1/(18m)}$ (which occurs with probability at least $1 - k_1\exp(-k_2n^{1/4})$ by Proposition \ref{prop:unifconv}).  Thus, we have that $\sup_{\theta\in\Theta}|R_n(\theta)-\overline{R}_n(\theta)| \stackrel{p}{\longrightarrow} 0$.  Now consider any $\hat{\theta}_n\in\arg\min\{\overline{R}_n(\theta)\ |\ \theta\in\Theta\}$, and note that the estimate $\hat{\theta}_n$ returned by the semiparametric algorithm satisfies this property by construction.  By definition we have $\overline{R}_n(\hat{\theta}_n) \leq \overline{R}_n(\theta_0)$, which can be rewritten as
\begin{equation}
R_n(\hat{\theta}_n) + \overline{R}_n(\hat{\theta}_n) - R_n(\hat{\theta}_n) \leq R_n(\theta_0) + \overline{R}_n(\theta_0) - R_n(\theta_0).
\end{equation}
Thus, we have
\begin{equation}
R_n(\hat{\theta}_n) \leq R_n(\theta_0) + |\overline{R}_n(\hat{\theta}_n) - R_n(\hat{\theta}_n)| + |\overline{R}_n(\theta_0) - R_n(\theta_0)|.
\end{equation}
We have thus shown all the conditions required to apply Theorem 5.14 of \citep{van2000}, which gives $\hat{\theta}_n \stackrel{p}{\longrightarrow} \theta_0$.  Now let $\overline{\theta}_n \in \arg\min\{Q_n(\theta)\ |\ \theta\in\Theta\}$.  By Theorem \ref{theorem:estconst}, we have $\overline{\theta}_n\stackrel{p}{\longrightarrow} \theta_0$.  This means that $|\overline{\theta}_n-\hat{\theta}_n|\stackrel{p}{\longrightarrow} 0$.
\Halmos\endproof
Theorem 7 states that the semiparametric algorithm produces estimates that are statistically consistent under the appropriate conditions. In the next section, we present several numerical experiments which validate our theoretical results as well as the performance of the enumeration and semiparametric algorithms.

\section{Numerical Experiments}
\label{section:data}
We present numerical results that demonstrate the statistical consistency of our algorithms for inverse optimization with noisy data, and the results show our algorithms perform competitively against \textsf{KKA} \citep{keshavarz2011} and \textsf{VIA} \citep{bertsimas2013}.  We begin by conducting two types of tests using synthetic data.  The first type is where the model is kept fixed and the number of data points increases, and the purpose is to demonstrate either estimation consistency or risk consistency of our algorithms.  The second type is where the number of data points is kept fixed and the number of the parameters in the model increases, and the purpose is to demonstrate the feasibility of using our algorithms on large-scale problems. We then apply our framework to a real data set, where we estimate a utility function that describes the tradeoff made between occupant comfort and energy consumption when setting a thermostat temperature setpoint for air-conditioning.
	
\subsection{Synthetic Data and Enumeration Algorithm}\label{sec:synthetic}

In the first experiments, we generate data using a given \textsf{FOP} and then use the same set of equations in \textsf{SAA-IOP}.  In other words, the first set of experiments are situations where the model whose parameters are being identified exactly match the model that generates the data.  As a result, this setting consists of situations where \textbf{IC} is satisfied.  The first example is where: (i) \textsf{FOP-A} is $\min \{(\theta+u)\cdot x\ |\ x\in[-1,1]\}$, (ii) $u$ has a uniform distribution with support $[-1,1]$, (iii) the measurement noise $w$ has a normal distribution with zero mean and unit variance, (iv) the data is generated with $\theta_0 = 1$, and (v) the enumeration algorithm (i.e., Algorithm \ref{algorithm:enum}) was applied with $\epsilon = 0.001$, $\delta = 0.01$, and $\Theta = [-1,1]$.  The second example is where: (i) \textsf{FOP-B} is $\min \{x^2-(\theta+u)\cdot x\ |\ x\in[0,1]\}$, (ii) $u$ has a uniform distribution with support $[0,2]$, (iii) the measurement noise $w$ has a normal distribution with zero mean and unit variance, (iv) the data is generated with $\theta_0 = \frac{1}{2}$, and (v) the enumeration algorithm (i.e., Algorithm \ref{algorithm:enum}) was applied with $\epsilon = 0$, $\delta = 0.01$, and $\Theta = [0,2]$.  

\begin{figure}
\centering
 \subfloat[\textsf{ENA}]{\includegraphics{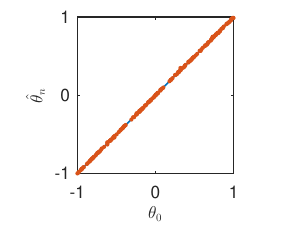}}
 \subfloat[\textsf{KKA}]{\includegraphics{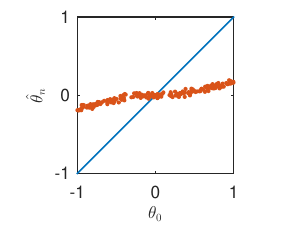}}
 \subfloat[\textsf{VIA}]{\includegraphics{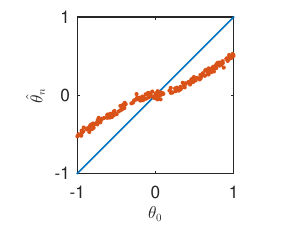}}\\
\caption{Scatter plot comparing estimated parameter $\hat{\theta}_n$ versus true parameter $\theta_0$ as computed by ENA, KKA and VIA algorithms at $n = 1,000$, when the data and model are both \textsf{FOP-A}.}
\label{fig:lplex}
\end{figure}

\begin{figure}
\centering
 \subfloat[\textsf{ENA}]{\includegraphics[trim=0.20in 0 0.23in 0, clip]{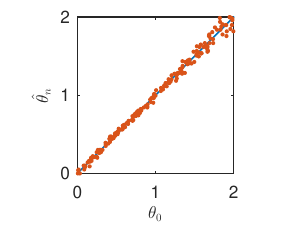}}
 \subfloat[\textsf{KKA}]{\includegraphics[trim=0.20in 0 0.23in 0, clip]{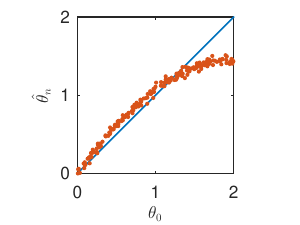}}
 \subfloat[\textsf{VIA}]{\includegraphics[trim=0.20in 0 0.23in 0, clip]{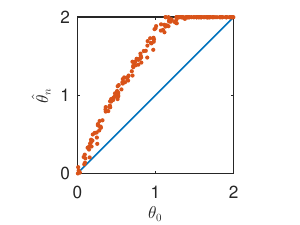}}
 \subfloat[\textsf{SPA}]{\includegraphics[trim=0.20in 0 0.23in 0, clip]{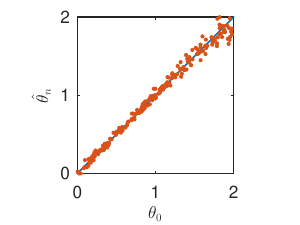}}\\
\caption{Scatter plot comparing estimated parameter $\hat{\theta}_n$ versus true parameter $\theta_0$ as computed by ENA, KKA, VIA and SPA algorithms at $n = 10,000$ when the data and model are both \textsf{FOP-B}.}
\label{fig:qplex}
\end{figure}

The results averaged over 100 repetitions of sampling $n \in \{10, 30, 50, 100, 300, 500, 1000\}$ data points and then estimating the parameter $\theta$ are summarized in Table \ref{table:syn1}.  We label the enumeration algorithm (i.e., Algorithm \ref{algorithm:enum}) as \textsf{ENA} in the table.  These results display estimation consistency of the enumeration algorithm since estimation error is decreasing to zero.  To further illustrate estimation consistency, we conducted an experiment with the two examples above where the data was generated with a $\theta_0$ that was randomly chosen from a uniform distribution with support $[-1,1]$ and $[0,2]$ for the first and second examples, respectively.  A plot comparing the estimates $\hat{\theta}_n$ to the true parameter $\theta_0$ for the first situation when $n = 1,000$ is shown in Figure \ref{fig:lplex}, and a plot comparing the estimates $\hat{\theta}_n$ to the true parameter $\theta_0$ for the second situation when $n = 10,000$ is shown in Figure \ref{fig:qplex}.  Consistent estimates should line up along the diagonal, and hence these plots demonstrate the estimation consistency (inconsistency) of the enumeration algorithm (\textsf{KKA} and \textsf{VIA}).  Recall from the discussion in Section 2 that \textsf{KKA} and \textsf{VIA} are inconsistent because they minimize an incorrect measure of error, and this discrepancy is most significant for points where the optimal solution of \textsf{FOP} lies on the boundary of the feasible set.  \textsf{KKA} and \textsf{VIA} perform more poorly for \textsf{FOP}-\textsf{A} than for \textsf{FOP}-\textsf{B} because \textsf{FOP}-\textsf{A} is a linear program, which has almost all of its optimal solutions on the boundary of the feasible set, whereas \textsf{FOP}-\textsf{B} is a quadratic program, which has more optimal solutions within the strict interior of the feasible set.

\begin{table}
\centering
\caption{\label{table:syn1} Estimation error $|\hat{\theta}_n - \theta_0|$ of enumeration algorithm (ENA) and benchmark algorithms (KKA and VIA) on two synthetic instances ($n$ increasing, $p = 1$).}
\begin{tabular}{l|cccccccc}
\hline
&$n$ & 10 & 30 & 50 & 100 & 300 & 500 & 1000\\
\hline
\multirow{3}{*}{\begin{tabular}{ll}Data:& \textsf{FOP-A}\\ Model:& \textsf{FOP-A}\end{tabular}}&\textsf{ENA} & 0.2616 & 0.0926 & 0.0380 & 0.0211 & 0.0055 & 0.0030 & 0.0009 \\
&\textsf{KKA} & 0.8686 & 0.8293 & 0.8182 & 0.8257 & 0.8130 & 0.8231 & 0.8170 \\
&\textsf{VIA} & 0.5552 & 0.4976 & 0.4829 & 0.4887 & 0.4807 & 0.4846 & 0.4780 \\
\hline
\multirow{3}{*}{\begin{tabular}{ll}Data:& \textsf{FOP-B}\\ Model:& \textsf{FOP-B}\end{tabular}}&\textsf{ENA} & 0.4577 & 0.2481 & 0.1510 & 0.0501 & 0.0222 & 0.0123 & 0.0063 \\
&\textsf{KKA} & 0.5065 & 0.2281 & 0.1595 & 0.0751 & 0.0398 & 0.0342 & 0.0238 \\
&\textsf{VIA} & 0.9488 & 0.7051 & 0.6344 & 0.4284 & 0.3145 & 0.3810 & 0.2962\\
\hline
\end{tabular}
\end{table}

In the second set of experiments, we generate data using a given model that is different than the \textsf{FOP} used to formulate \textsf{SAA-IOP}.  In other words, this set of experiments are situations where the model whose parameters are being identified does not match the model that generates the data.  As a result, this setting consists of situations where \textbf{IC} is \emph{not} satisfied.  The first example is where: (i) the data is generated by \textsf{FOP-C} which is $\min \{\frac{3}{2}\cdot x^2-(1+u)\cdot x\ |\ x\in[0,1]\}$, (ii) the model estimated by \textsf{IOP-SAA} is \textsf{FOP-B}, (iii) $u$ has a uniform distribution with support $[0,5]$, (iv) the measurement noise $w$ has a normal distribution with zero mean and unit variance, and (v) the enumeration algorithm (i.e., Algorithm \ref{algorithm:enum}) was applied with $\epsilon = 0$, $\delta = 0.01$, and $\Theta = [0,2]$.  The second example is where: (i) the data is generated by the statistical model \textsf{SQR-1} given by $y_i = \min\{\max\{\sqrt{u_i},0\},1\} + w_i$, (ii) the model estimated by \textsf{IOP-SAA} is \textsf{FOP-B}, (iii) $u$ has a uniform distribution with support $[0,5]$, (iv) the measurement noise $w$ has a normal distribution with zero mean and unit variance, and (v) the enumeration algorithm (i.e., Algorithm \ref{algorithm:enum}) was applied with $\epsilon = 0$, $\delta = 0.01$, and $\Theta = [0,2]$.  

The results averaged over 100 repetitions of sampling $n \in \{10, 30, 50, 100, 300, 500, 1000\}$ data points and then estimating the parameter $\theta$ are summarized in Table \ref{table:syn2}, and these results are normalized by subtracting $\text{var}(w)$.  The reason for this normalization is that the prediction error $\mathbb{E}((y-\xi(u))^2))$ of the prediction $\xi(u)$ of the true model (either \textsf{FOP-C} or \textsf{SQR-M}, respectively) is $\text{var}(w)$ because $y = \xi(u) + w$ here.  The enumeration algorithm has lower prediction error because it is risk consistent, whereas \textsf{KKA} and \textsf{VIA} are not risk consistent.

\begin{table}
\centering
\caption{\label{table:syn2}Normalized prediction error $Q(\hat{\theta}_n)-\text{var}(w)$ of enumeration algorithm (ENA) and benchmark algorithms (KKA and VIA) on two synthetic instances ($n$ increasing, $p = 1$).}\begin{tabular}{l|cccccccc}
\hline
&$n$ & 10 & 30 & 50 & 100 & 300 & 500 & 1000\\
\hline
\multirow{3}{*}{\begin{tabular}{ll}Data:& \textsf{FOP-C}\\ Model:& \textsf{FOP-B}\end{tabular}}&\textsf{ENA} & 0.0216 & 0.0184 & 0.0162 & 0.0150 & 0.0065 & 0.0046 & 0.0017 \\
&\textsf{KKA} & 0.0168 & 0.0124 & 0.0128 & 0.0151 & 0.0150 & 0.0150 & 0.0132\\
&\textsf{VIA} & 0.0249 & 0.0185 & 0.0196 & 0.0149 & 0.0089 & 0.0072 & 0.0042 \\
\hline
\multirow{3}{*}{\begin{tabular}{ll}Data:& \textsf{SQR-1}\\ Model:& \textsf{FOP-B}\end{tabular}}&\textsf{ENA} & 0.0294 & 0.0217 & 0.0152 & 0.0110 & 0.0073&0.0041 &0.0024 \\
&\textsf{KKA} & 0.0394 & 0.0389 & 0.0398 & 0.0440 & 0.0504 & 0.0525 & 0.0518 \\
&\textsf{VIA} & 0.0343 & 0.0287 & 0.0243 & 0.0187 & 0.0122 & 0.0084 & 0.0072 \\
\hline
\end{tabular}
\end{table}

\subsection{Synthetic Data and Semiparametric Algorithm}\label{sec:spac}
We now examine the performance of the semiparametric algorithm (Algorithm 2) in four sets of experiments. In the first set of experiments, we generate data using a given \textsf{FOP} and then use the same equations in \textsf{SAA-IOP}.  These experiments are situations where the model whose parameters are being identified exactly matches the model that generates the data.  As a result, this setting consists of situations where \textbf{IC} is satisfied. We consider three different formulations for \textsf{FOP}.  The first example is where: (i) \textsf{FOP-D} is $\min \{x'x-(\theta+u)'x\ |\ x\in[0,1]^p\}$, (ii) $u$ has a uniform distribution with support $[0,2]^p$, (iii) the measurement noise $w$ has a jointly Gaussian distribution with zero mean and identity covariance, (iv) the data is generated with $p = 10$ and $\theta_0\in\mathbb{R}^p$ such that $ \theta_{0k} = \frac{1}{2}$ for all $k\in[p]$, and (v) the semiparametric algorithm (i.e., Algorithm \ref{algorithm:sp}) was applied with $\gamma,\sigma$ chosen using cross-validation \citep{hastie2009} and $\Theta = [0,2]$.  The second example is where: (i) \textsf{FOP-E} is 
\begin{equation}
\min \left\{ -\textstyle\sum_{k=1}^{p} \theta_k\cdot\log( x_k +  u_k) - \log( x_{p+1} +  u_{p+1})  ~\big|~  x_k \geq 0,\ \textstyle\sum_{k=1}^{p+1} x_k = 1 \right\},
\end{equation}
(ii) $u$ has a uniform distribution with support $[1,2]^{p+1}$, (iii) the measurement noise $w$ has a jointly Gaussian distribution with zero mean and identity covariance, (iv) the data is generated with $p = 10$ and $\theta_0\in\mathbb{R}^p$ such that $ \theta_{0k} = 1$ for all $k\in[p]$, and (v) a modified version of the seimparametric algorithm (i.e., Algorithm \ref{algorithm:sp}) was applied with $\gamma,\sigma$ chosen using cross-validation  and $\Theta = [\frac{1}{2},2]$. The modification to Algorithm \ref{algorithm:sp} is that we calculate $\tilde{x}_i = \min_x\{\|\overline{x}_i - x\|\ |\  x_k \geq 0\}$ and then compute $\hat{\theta}_n$ using $\mathsf{SP}$--$\mathsf{IOP}$--$\mathsf{RISK}$--$\mathsf{SAA}$ with the $\tilde{x}_i$ replacing the $\overline{x}_i$.  The $\tilde{x}_i$ are the projection of the $\overline{x}_i$ onto the nonnegative orthant, and it turns out this projection does not affect our theoretical results.  In particular, a short proof using the continuous mapping theorem \citep{van2000} and the boundedness of the feasible set in \textbf{R1} gives that $\max_{i\in[n]}\big\|\tilde{x}_i - \mathcal{S}(u_i,\theta_0)\big\|\stackrel{p}{\longrightarrow} 0$.  The projection is needed for this particular example because otherwise the inverse formulation would contain logarithms of negative numbers, which are complex-valued.  More generally, a projection of $\overline{x}_i$ onto the feasible set of \textsf{FOP} will not affect our theoretical results, and can be added as a step in our semiparametric algorithm.

\begin{figure}
\centering
 \subfloat[\textsf{ENA}]{\includegraphics[trim=0.20in 0 0.23in 0, clip]{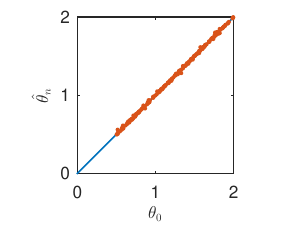}}
 \subfloat[\textsf{KKA}]{\includegraphics[trim=0.20in 0 0.23in 0, clip]{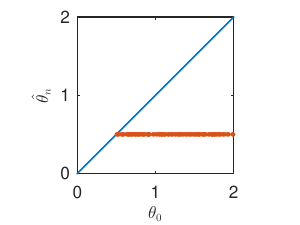}}
 \subfloat[\textsf{VIA}]{\includegraphics[trim=0.20in 0 0.23in 0, clip]{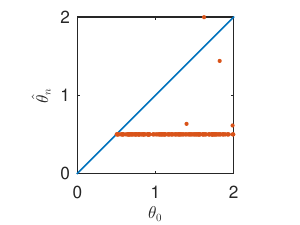}}
 \subfloat[\textsf{SPA}]{\includegraphics[trim=0.20in 0 0.23in 0, clip]{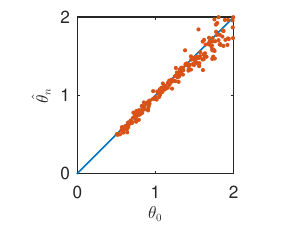}}\\
\caption{Scatter plot comparing estimated parameter $\hat{\theta}_n$ versus true parameter $\theta_0$ as computed by different algorithms at $n = 1,000$ when the data and model are both \textsf{FOP-E}.}
\label{fig:logplex}
\end{figure}

In Table \ref{table:syn3}, we present estimation results for the first and second examples, averaged over 100 repetitions for each value of $n \in \{10, 30, 50, 100, 300, 500, 1000\}$.  We label the semiparametric algorithm (i.e., Algorithm \ref{algorithm:sp}) as \textsf{SPA} in the table.  These results display estimation consistency of the semiparametric algorithm since it has lower estimation error as the data increases.  To further illustrate estimation consistency, we conducted an experiment with the two situations above where the data was generated with $p=1$ and a $\theta_0$ that was randomly chosen from a uniform distribution with support $[0,1]$ and $[\frac{1}{2},2]$ for the first and second situations, respectively.  A plot comparing the estimates $\hat{\theta}_n$ to the true parameter $\theta_0$ for the first situation when $n = 1,000$ is shown in Figure \ref{fig:qplex}, and a plot comparing the estimates $\hat{\theta}_n$ to the true parameter $\theta_0$ for the second situation when $n = 1,000$ is shown in Figure \ref{fig:logplex}.  Consistent estimates should line up along the diagonal, and hence these plots demonstrate the estimation consistency (inconsistency) of the semiparametric algorithm (\textsf{KKA} and \textsf{VIA}).  It is worth comparing the results of the semiparametric and enumeration algorithms.  As mentioned above, the semiparametric algorithm will generally have higher estimation error than the enumeration algorithm -- this can be observed in these plots because the semiparametric algorithm estimates have a larger variation about the diagonal than the estimates of the enumeration algorithm. 

\begin{table}
\centering
\caption{\label{table:syn3}Estimation error $\|\hat{\theta}_n - \theta_0\|$ of semiparametric algorithm (SPA) and benchmark algorithms (KKA and VIA) on two synthetic instances ($n$ increasing, $p = 10$).}
\begin{tabular}{l|cccccccc}
\hline
&$n$ & 10 & 30 & 50 & 100 & 300 & 500 & 1000\\
\hline
\multirow{3}{*}{\begin{tabular}{ll}Data:& \textsf{FOP-D}\\ Model:& \textsf{FOP-D}\end{tabular}}&\textsf{SPA} & 2.4618 & 1.7025 & 1.2543 & 0.8535 & 0.4754 & 0.3750 & 0.2573 \\
&\textsf{KKA} & 2.2569 & 1.5513 & 1.2229 & 0.9281 & 0.6107 & 0.5435 & 0.4447 \\
&\textsf{VIA} & 3.3829 & 3.2603 & 3.1937 & 3.1501 & 3.0292 & 3.0324 & 2.9208 \\
\hline
\multirow{3}{*}{\begin{tabular}{ll}Data:& \textsf{FOP-E}\\ Model:& \textsf{FOP-E}\end{tabular}}&\textsf{SPA} & 0.9189 & 0.7982 & 0.7500 & 0.7487 & 0.6639 & 0.6070 & 0.5783 \\
&\textsf{KKA} & 1.6687 & 1.5850 & 1.5813 & 1.5865 & 1.5828 & 1.5806 & 1.5811 \\
&\textsf{VIA} & 1.9299 & 1.6781 & 1.6826 & 1.6132 & 1.6001 & 1.5973 & 1.5843\\
\hline
\end{tabular}
\end{table}

In the second set of experiments, we generate data using a given model that is different than the \textsf{FOP} used to formulate \textsf{SAA-IOP}.  In other words, this set of experiments are situations where the model whose parameters are being identified does not match the model that generates the data.  As a result, this setting consists of situations where \textbf{IC} is \emph{not} satisfied.  The first setting is where: (i) the data is generated by \textsf{FOP-C} which is $\min \{\frac{3}{2}\cdot x'x-(1+u)'x\ |\ x\in[0,1]^{10}\}$, (ii) the model estimated by \textsf{IOP-SAA} is \textsf{FOP-D} with $p = 10$, (iii) $u$ has a uniform distribution with support $[0,5]^{10}$, (iv) the measurement noise $w$ has a jointly Gaussian distribution with zero mean and identity covariance, and (v) the semiparametric algorithm (i.e., Algorithm \ref{algorithm:sp}) was applied with $\gamma,\sigma$ chosen using cross-validation \citep{hastie2009} and $\Theta = [0,2]$.  The second setting is where: (i) the data is generated by the statistical model \textsf{SQR-P} given by $y_i = \min\{\max\{\sqrt{u_i},0\},1\} + w_i$, (ii) the model estimated by \textsf{IOP-SAA} is \textsf{FOP-D} with $p = 10$, (iii) $u$ has a uniform distribution with support $[0,5]^{10}$, (iv) the measurement noise $w$ has a jointly Gaussian distribution with zero mean and identity covariance, and (v) the semiparametric algorithm (i.e., Algorithm \ref{algorithm:sp}) was applied with $\gamma,\sigma$ chosen using cross-validation \citep{hastie2009} and $\Theta = [0,2]$. The results averaged over 100 repetitions of sampling $n \in \{10, 30, 50, 100, 300, 500, 1000\}$ data points and then estimating the parameter $\theta$ are summarized in Table \ref{table:syn4}, and these results are normalized by subtracting $\mathbb{E}(w'w)$.  The reason for this normalization is that the prediction error $\mathbb{E}(\|y-\xi(u)\|^2))$ of the prediction $\xi(u)$ of the true model (either \textsf{FOP-C} or \textsf{SQR-M}, respectively) is $\mathbb{E}(w'w)$ because $y = \xi(u) + w$ here.  The enumeration algorithm has lower prediction error because it is risk consistent, whereas \textsf{KKA} and \textsf{VIA} are not risk consistent.

\begin{table}
\centering
\caption{\label{table:syn4}Normalized prediction error $Q(\hat{\theta}_n)-\mathbb{E}(w'w)$ of semiparametric algorithm (SPA) and benchmark algorithms (KKA and VIA) on two synthetic instances ($n$ increasing, $p = 10$).}
\begin{tabular}{l|cccccccc}
\hline
&$n$ & 10 & 30 & 50 & 100 & 300 & 500 & 1000\\
\hline
\multirow{3}{*}{\begin{tabular}{ll}Data:& \textsf{FOP-C}\\ Model:& \textsf{FOP-D}\end{tabular}}&\textsf{SPA} & 0.2319 & 0.1972 & 0.1744 & 0.1501 & 0.1029 & 0.0844 & 0.0529 \\
&\textsf{KKA} & 0.1584 & 0.1308 & 0.1314 & 0.1349 & 0.1452 & 0.1497 & 0.1481\\
&\textsf{VIA} & 0.3438 & 0.3407 & 0.3360 & 0.3205 & 0.2950 & 0.2816 & 0.2811 \\
\hline
\multirow{3}{*}{\begin{tabular}{ll}Data:& \textsf{SQR-M}\\ Model:& \textsf{FOP-D}\end{tabular}}&\textsf{SPA} & 0.4180 & 0.3497 & 0.3195 & 0.2470 & 0.1572 & 0.0998 &0.0658 \\
&\textsf{KKA} & 0.3645 & 0.3885 & 0.3987 & 0.4537 & 0.5115 & 0.5114 & 0.5214 \\
&\textsf{VIA} & 0.3468 & 0.2784 & 0.2737 & 0.2524 & 0.2405 & 0.2458 & 0.2599 \\
\hline
\end{tabular}
\end{table}

In the third set of experiments, we generate data using the previous four settings.  The difference in this set of experiments is that we fix $n=1,000$ and vary $p\in\{1,3,5,10,30\}$.  The results when the data/model are given by \textsf{FOP-D}/\textsf{FOP-D} and \textsf{FOP-E}/\textsf{FOP-E}, averaged over 100 repetitions and then estimating the parameter $\theta$, are summarized in Table \ref{table:syn5}.  These results show that the semiparametric algorithm has lower estimation error than \textsf{KKA} and \textsf{VIA} on these examples. The results when the data/model are given by \textsf{FOP-C}/\textsf{FOP-B} and \textsf{SQR-M}/\textsf{FOP-B}, averaged over 100 repetitions and then estimating the parameter $\theta$, are summarized in Table \ref{table:syn6}.  These results show that the semiparametric algorithm has lower prediction error than \textsf{KKA} and \textsf{VIA} on these examples.

\begin{table}
\centering
\caption{\label{table:syn5}Estimation error $\|\hat{\theta}_n-\theta_0\|$ of semiparametric algorithm (SPA) and benchmark algorithms (KKA and VIA) on two synthetic instances ($n=1,000$, $p$ increasing).}
\begin{tabular}{l|cccccccc}
\hline
&$p$ & 1 & 3 & 5 & 10 & 30 \\
\hline
\multirow{3}{*}{\begin{tabular}{ll}Data:& \textsf{FOP-D}\\ Model:& \textsf{FOP-D}\end{tabular}}&\textsf{SPA} & 0.0601 & 0.1464 & 0.1907 & 0.2794 & 0.4701 \\
&\textsf{KKA} & 0.1178 & 0.2349 & 0.3038 & 0.4619 & 0.7978 \\
&\textsf{VIA} & 0.4943 & 1.2254 & 1.8099 & 2.9522 & 5.7737 \\
\hline
\multirow{3}{*}{\begin{tabular}{ll}Data:& \textsf{FOP-E}\\ Model:& \textsf{FOP-E}\end{tabular}}&\textsf{SPA} & 0.0251 & 0.1258 & 0.2571 & 0.5890 & 0.5576 \\
&\textsf{KKA} & 0.5000 & 0.8660 & 1.1174 & 1.5804 & 2.7377 \\
&\textsf{VIA} & 0.5000 & 0.8691 & 1.1231 & 1.5966 & 2.7628\\
\hline
\end{tabular}
\end{table}

\begin{table}
\centering
\caption{\label{table:syn6}Normalized prediction error $Q(\hat{\theta}_n)-\mathbb{E}(w'w)$ of semiparametric algorithm (SPA) and benchmark algorithms (KKA and VIA) on two synthetic instances ($n=1,000$, $p$ increasing).}
\begin{tabular}{l|cccccccc}
\hline
&$p$ & 1 & 3 & 5 & 10 & 30 \\
\hline
\multirow{3}{*}{\begin{tabular}{ll}Data:& \textsf{FOP-C}\\ Model:& \textsf{FOP-D}\end{tabular}}&\textsf{SPA} & 0.0064 & 0.0171 & 0.0403 & 0.0628 & 0.2048 \\
&\textsf{KKA} & 0.0538 & 0.1553 & 0.2619 & 0.5252 & 1.5712\\
&\textsf{VIA} & 0.0078 & 0.0175 & 0.0745 & 0.2602 & 0.9654 \\
\hline
\multirow{3}{*}{\begin{tabular}{ll}Data:& \textsf{SQR-M}\\ Model:& \textsf{FOP-D}\end{tabular}}&\textsf{SPA} & 0.0056 & 0.0194 & 0.0319 & 0.0606 & 0.1568 \\
&\textsf{KKA} & 0.0148 & 0.0471 & 0.0761 & 0.1523 & 0.4394 \\
&\textsf{VIA} & 0.0055 & 0.0273 & 0.0821 & 0.2848 & 1.2896 \\
\hline
\end{tabular}
\end{table}
\subsection{High-Dimensional Nonlinear Forward Problem with Stochastic Constraints}
We now consider a setting where $\textsf{FOP}$ is high dimensional, contains a logarithmic objective, and has an exponential stochastic constraint (i.e. the constraint depends on $u$). Specifically, we consider the following setting: (i) \textsf{FOP-F} is 
\begin{equation}
\min \left\{ -\textstyle\sum_{k=1}^{p} \theta_k\cdot  u^{(1)}_k\cdot\log( x_k )  ~\big|~ \frac{1}{p}\sum_{k=1}^p e^{x_k + u^{(1)}_k} - u^{(2)}_k \le 0, \;  x_k \geq 0 \right\},
\end{equation}
(ii) $u^{(1)}$ has a uniform distribution with support $[1,2]^{p}$ and $u^{(2)}$ has a uniform distribution with support $[50,100]^{p}$, (iii) the measurement noise $w$ has a jointly Gaussian distribution with zero mean and identity covariance, (iv) the data is generated with $\theta_0\in\mathbb{R}^p_+$ such that $ \sum_{k=1}^p  \theta_{0k} = p$, and (v) a modified version of the seimparametric algorithm (i.e., Algorithm \ref{algorithm:sp}) is applied where $\Theta = \{\theta \in \mathbb{R}_+^p  | \sum_{k=1}^p  \theta_{k} = p \}$ and $\gamma,\sigma$ is selected using cross-validation. We set $n = 1,000$ and repeat the sampling and estimation procedure 100 times for each value of $p \in \{5,10,20,50,100\}$. The average estimation and prediction errors are summarized in Tables \ref{table:syn7} and \ref{table:syn8}, respectively, which show that the semiparametric algorithm is competitive with existing methods in this setting as well. Note that the magnitude of the errors is expected to increase with $p$, since we do not normalize the error for the number of parameters being estimated.

\begin{table}
\centering
\caption{\label{table:syn7}Estimation error $\|\hat{\theta}_n-\theta_0\|$ of semiparametric algorithm (SPA) and benchmark algorithms (KKA and VIA) on synthetic instance ($n = 1,000$, $p$ increasing).}
\begin{tabular}{l|cccccccc}
\hline
&$p$ & 5 & 10 & 20 & 50 & 100 \\
\hline
\multirow{3}{*}{\begin{tabular}{ll}Data:& \textsf{FOP-F}\\ Model:& \textsf{FOP-F}\end{tabular}}&\textsf{SPA} & 0.5535 &   0.8530 &   1.1522 &    2.0020 &   2.7205  \\
&\textsf{KKA} &2.1753   & 4.6199    &8.4599  & 12.4102   & 17.8112  \\
&\textsf{VIA} & 1.1825  &  1.8689   & 3.7320 &   5.9003   & 8.1874  \\
\hline
\end{tabular}
\end{table}

\begin{table}
\centering
\caption{\label{table:syn8}Normalized prediction error $Q(\hat{\theta}_n)-\mathbb{E}(w'w)$ of semiparametric algorithm (SPA) and benchmark algorithms (KKA and VIA) on synthetic instance ($n = 1,000$, $p$ increasing).}
\begin{tabular}{l|rrrrrrrr}
\hline
& \multicolumn{1}{c}{$p$} & \multicolumn{1}{c}{5} & \multicolumn{1}{c}{10} & \multicolumn{1}{c}{20} & \multicolumn{1}{c}{50} & \multicolumn{1}{c}{100} \\
\hline
\multirow{3}{*}{\begin{tabular}{ll}Data:& \textsf{FOP-F}\\ Model:& \textsf{FOP-F}\end{tabular}}&\textsf{SPA}& 0.2539 &   0.5117 &   0.9307 &   2.9423 &   5.8644  \\
&\textsf{KKA} & 8.2329 &  22.5149 &  60.4496  & 145.2210 & 302.1124 \\
&\textsf{VIA} & 2.3475  &  3.7495 &  10.8325  & 26.5713  & 48.3217  \\
\hline
\end{tabular}
\end{table}

\FloatBarrier

\subsection{Empirical Data: Estimating an Energy-Comfort Utility Function}\label{sec:vehicle}
We next apply our inverse optimization framework to the problem of estimating a utility function that describes the tradeoff made between occupant comfort and the amount of energy consumption when setting a thermostat temperature setpoint for air-conditioning.  The data we use is collected from Sutardja Dai Hall on the Berkeley campus, which was used as part of the BRITE-S testbed in our past experiments \citep{aswani2012_ifac,aswani2012_valuetools,aswani2012_acc2} concerning robust learning-based optimization \citep{aswani2013_automatica} of heating, ventilation, and air-conditioning (HVAC) systems.  Specifically, this building is equipped with a commercial web application \citep{comfy2016} that allows occupants to change the thermostat temperature setpoints in real-time, and so the setpoints are changed throughout the year by occupants in response to factors like the outside weather.

When a room is being cooled, a lower temperature setpoint requires increased energy consumption since the air-conditioner must provide more cold air; however, the purpose of air-conditioning is to improve comfort by lowering the room temperature.  And so individuals must tradeoff comfort and energy consumption when choosing the setpoint.  A simplified utility function model (expressed as minimization of the negative of the utility function) that captures this tradeoff is \textsf{FOP-S}:
\begin{equation}
\min_x \big\{ \theta_1\cdot(x-76)^2 + (x- \theta_2-u)^2\ \big|\ x\in[70,76]\big\},
\end{equation}
where $x \in \mathbb{R}$ is the thermostat temperature setpoint in units of degrees Fahrenheit ($^\circ$F), and $u\in\mathbb{R}$ is the current outside temperaure in degrees Fahrenheit ($^\circ$F).  The term $(x- \theta_2-u)^2$ indicates a preference for a temperature setpoint that is a fixed amount $ \theta_2$ above the outside temperature $u$ (i.e., the preferred temperature is $ \theta_2 + u$), and the reason for this term is that individuals prefer a higher indoor temperature as the outside temperature increases \citep{ashrae2013}.  The term $ \theta_1\cdot(x-76)^2$ indicates a preference for a higher setpoint because of energy considerations, and the number 76 is used because 76$^\circ$F--78$^\circ$F is a relatively high setpoint temperature that is often recommended for saving energy.  The parameter $ \theta_1$ quantifies the tradeoff between the preference for a higher setpoint to save energy versus the desired indoor temperature $ \theta_2 + u$.  Lastly, the constraints $x\in[70,76]$ indicate observed setpoint limits.

The results averaged over 100 repetitions of sampling $n \in \{10, 30, 50, 100, 300, 500, 1000\}$ data points and then estimating the parameters $\theta$ are summarized in Table \ref{table:sdh}.  The data set (which we label \textsf{SDH-E} in the table) used consists of outside temperature measurements (i.e., the $u$ variable) and the chosen temperature set point (i.e., the $x$ variable) of a single thermostat in Sutardja Dai Hall.  In each repetition, the full data set was randomly split into a 1,000 point \emph{training} data set and a 14,500 point \emph{testing} data set.  The $n$ data points were randomly chosen from the training data set, and the prediction error of the estimated parameters were computed using the testing data set.  To evaluate the statistical significance of the computed results, a bootstrap hypothesis test \citep{efron1994} was conducted.  The computed $p$-value was less than 0.01, which indicates that the improved performance of the enumeration algorithm is statistically significant.

\begin{table}
\centering
\caption{\label{table:sdh}Prediction error $Q(\hat{\theta}_n)$ of enumeration algorithm (ENA) and benchmark algorithms (KKA and VIA) on temperature preference dataset. ($n$ increasing, $p$ fixed).}
\begin{tabular}{l|cccccccc}
\hline
&$n$ & 10 & 30 & 50 & 100 & 300 & 500 & 1000\\
\hline
\multirow{3}{*}{\begin{tabular}{ll}Data:& \textsf{SDH-E}\\ Model:& \textsf{FOP-S}\end{tabular}}&\textsf{ENA} & 1.3656 & 1.3308 & 1.3255 & 1.3169 & 1.3112 & 1.3099 & 1.3090 \\
&\textsf{KKA} & 2.2439 & 2.2528 & 2.2508 & 2.2351 & 2.2225 & 2.2220 & 2.2200\\
&\textsf{VIA} & 2.2975 & 2.2538 & 2.2472 & 2.2277 & 2.2163 & 2.2166 & 2.2138 \\
\hline
\end{tabular}
\end{table}

\section{Conclusion}
We developed and analyzed a formulation for inverse optimization in the setting where noisy measurements of the optimal points of a convex optimization problem are available.  Our approach requires solving a bilevel program, and we defined a new duality-based reformulation to convert this bilevel program into a single level program. We showed that our formulation as a bilevel program leads to statistical consistency, in contrast to existing heuristics.  Although our formulation is NP-hard to solve, we provided two numerical algorithms that maintain the statistical consistency of our formulation. Finally, we demonstrated that our approach improves upon existing methods for inverse optimization through a series of numerical experiments using both synthetic and empirical data.

\ACKNOWLEDGMENT{The authors gratefully acknowledge the support of NSF Award CMMI-1450963 and an NSERC Postgraduate Scholarship.\\}

\noindent{\bf \large Authors}\\

\noindent{\bf Anil Aswani} is an Assistant Professor in the Department of Industrial Engineering and Operations Research at the University of California, Berkeley.  His research interests include data-driven decision making, with particular emphasis on addressing inefficiencies and inequities in health systems and physical infrastructure..

{\bf Zuo-Jun (Max) Shen} is a Chancellor's Professor in the Department of Industrial Engineering and Operations Research and the Department of Civil and Environmental Engineering at UC Berkeley. He is also an honorary professor at Tsinghua University. He has been active in the following research areas: integrated supply chain design and management, design and analysis of optimization algorithms, energy system and transportation system planning and optimization.

{\bf Auyon Siddiq} is a Ph.D. candidate in the Department of Industrial Engineering and Operations Research at the University of California, Berkeley. His current research interests are in data-driven optimization, healthcare operations and incentive design.\\



\bibliographystyle{ormsv080} 
\bibliography{articlesIO} 

\begin{thebibliography}{59}
\expandafter\ifx\csname natexlab\endcsname\relax\def\natexlab#1{#1}\fi
\expandafter\ifx\csname url\endcsname\relax
  \def\url#1{{\tt #1}}\fi
\expandafter\ifx\csname urlprefix\endcsname\relax\def\urlprefix{URL }\fi
\expandafter\ifx\csname urlstyle\endcsname\relax
  \expandafter\ifx\csname doi\endcsname\relax
  \def\doi#1{doi:\discretionary{}{}{}#1}\fi \else
  \expandafter\ifx\csname doi\endcsname\relax
  \def\doi{doi:\discretionary{}{}{}\begingroup \urlstyle{rm}\Url}\fi \fi

\bibitem[{Aalami et~al.(2010)Aalami, Moghaddam, and Yousefi}]{aalami2010demand}
Aalami, HA, M~Parsa Moghaddam, GR~Yousefi. 2010.
\newblock Demand response modeling considering interruptible/curtailable loads
  and capacity market programs.
\newblock {\it Applied Energy\/} {\bf 87}(1) 243--250.

\bibitem[{Ahuja and Orlin(2001)}]{ahuja2001}
Ahuja, Ravindra~K, James~B Orlin. 2001.
\newblock Inverse optimization.
\newblock {\it Operations Research\/} {\bf 49}(5) 771--783.

\bibitem[{ASHRAE(2013)}]{ashrae2013}
ASHRAE. 2013.
\newblock {\it ANSI/ASHRAE Standard 55-2013: Thermal Environmental Conditions
  for Human Occupancy\/}.
\newblock ASHRAE.

\bibitem[{Aswani(2015)}]{aswani2015}
Aswani, A. 2015.
\newblock Low-rank approximation and completion of positive tensors
  \urlprefix\url{http://arxiv.org/abs/1412.0620}.
\newblock ArXiv:1412.0620.

\bibitem[{Aswani et~al.(2013)Aswani, Gonzalez, Sastry, and
  Tomlin}]{aswani2013_automatica}
Aswani, A., H.~Gonzalez, S.~Sastry, C.~Tomlin. 2013.
\newblock Provably safe and robust learning--based model predictive control.
\newblock {\it Automatica\/} {\bf 49}(5) 1216--1226.

\bibitem[{Aswani et~al.(2016)Aswani, Kaminsky, Mintz, Flowers, and
  Fukuoka}]{aswani2016}
Aswani, A., P.~Kaminsky, Y.~Mintz, E.~Flowers, Y.~Fukuoka. 2016.
\newblock Predictive modeling of behavior in weight loss interventions
  Submitted.

\bibitem[{Aswani et~al.(2012{\natexlab{a}})Aswani, Master, Taneja, Krioukov,
  Culler, and Tomlin}]{aswani2012_ifac}
Aswani, A., N.~Master, J.~Taneja, A.~Krioukov, D.~Culler, C.~Tomlin.
  2012{\natexlab{a}}.
\newblock Energy--efficient building {HVAC} control using hybrid system
  {LBMPC}.
\newblock {\it IFAC Conference on Nonlinear Model Predictive Control\/}.

\bibitem[{Aswani et~al.(2012{\natexlab{b}})Aswani, Master, Taneja, Krioukov,
  Culler, and Tomlin}]{aswani2012_valuetools}
Aswani, A., N.~Master, J.~Taneja, A.~Krioukov, D.~Culler, C.~Tomlin.
  2012{\natexlab{b}}.
\newblock Quantitative methods for comparing different {HVAC} control schemes.
\newblock {\it International Conference on Performance Evaluation Methodologies
  and Tools\/}.

\bibitem[{Aswani et~al.(2012{\natexlab{c}})Aswani, Master, Taneja, Smith,
  Krioukov, Culler, and Tomlin}]{aswani2012_acc2}
Aswani, A., N.~Master, J.~Taneja, V.~Smith, A.~Krioukov, D.~Culler, C.~Tomlin.
  2012{\natexlab{c}}.
\newblock Identifying models of {HVAC} systems using semi-parametric
  regression.
\newblock {\it American Control Conference\/}.

\bibitem[{Aswani and Tomlin(2012)}]{aswani2012_allerton}
Aswani, A., C.~Tomlin. 2012.
\newblock Incentive design for efficient building quality of service.
\newblock {\it Allerton Conference on Communication, Control, and Computing\/}.
  90--97.

\bibitem[{Audet et~al.(1997)Audet, Hansen, Jaumard, and Savard}]{audet1997}
Audet, C., P.~Hansen, B.~Jaumard, G.~Savard. 1997.
\newblock Links between linear bilevel and mixed 0--1 programming problems.
\newblock {\it Journal of Optimization Theory and Applications\/} {\bf 93}(2)
  273--300.
\newblock \doi{10.1023/A:1022645805569}.
\newblock \urlprefix\url{http://dx.doi.org/10.1023/A:1022645805569}.

\bibitem[{Bajari et~al.(2007)Bajari, Benkard, and Levin}]{bajari2007}
Bajari, Patrick, C~Lanier Benkard, Jonathan Levin. 2007.
\newblock Estimating dynamic models of imperfect competition.
\newblock {\it Econometrica\/} {\bf 75}(5) 1331--1370.

\bibitem[{Bard and Moore(1990)}]{bard1990branch}
Bard, Jonathan~F, James~T Moore. 1990.
\newblock A branch and bound algorithm for the bilevel programming problem.
\newblock {\it SIAM Journal on Scientific and Statistical Computing\/} {\bf
  11}(2) 281--292.

\bibitem[{Bartlett and Mendelson(2002)}]{bartlett2002}
Bartlett, P., S.~Mendelson. 2002.
\newblock Rademacher and gaussian complexities: Risk bounds and structural
  results.
\newblock {\it J. Mach. Learn. Res.\/} .

\bibitem[{Beil and Wein(2003)}]{beil2003}
Beil, Damian~R, Lawrence~M Wein. 2003.
\newblock An inverse-optimization-based auction mechanism to support a
  multiattribute rfq process.
\newblock {\it Management Science\/} {\bf 49}(11) 1529--1545.

\bibitem[{Berge(1963)}]{berge1963}
Berge, Claude. 1963.
\newblock {\it Topological Spaces: including a treatment of multi-valued
  functions, vector spaces, and convexity\/}.
\newblock Courier Dover Publications.

\bibitem[{Bertsimas et~al.(2012)Bertsimas, Gupta, and
  Paschalidis}]{bertsimas2012}
Bertsimas, Dimitris, Vishal Gupta, Ioannis~Ch Paschalidis. 2012.
\newblock Inverse optimization: a new perspective on the black-litterman model.
\newblock {\it Operations research\/} {\bf 60}(6) 1389--1403.

\bibitem[{Bertsimas et~al.(2015)Bertsimas, Gupta, and
  Paschalidis}]{bertsimas2013}
Bertsimas, Dimitris, Vishal Gupta, Ioannis~Ch Paschalidis. 2015.
\newblock Data-driven estimation in equilibrium using inverse optimization.
\newblock {\it Mathematical Programming\/} {\bf 153}(2) 595--633.

\bibitem[{Bickel and Doksum(2006)}]{bickel2006}
Bickel, P., K.~Doksum. 2006.
\newblock {\it Mathematical Statistics: Basic Ideas And Selected Topics\/},
  vol.~1.
\newblock 2nd ed. Pearson Prentice Hall.

\bibitem[{Bonnans and Shapiro(2000)}]{bonnans2000}
Bonnans, J., A.~Shapiro. 2000.
\newblock {\it Perturbation Analysis of Optimization Problems\/}.
\newblock Springer.

\bibitem[{Bonnans and Ioffe(1995{\natexlab{a}})}]{bonnans1995}
Bonnans, J~Frederic, Alexander~D Ioffe. 1995{\natexlab{a}}.
\newblock Quadratic growth and stability in convex programming problems with
  multiple solutions.
\newblock {\it J. Convex Anal\/} {\bf 2}(1-2) 41--57.

\bibitem[{Bonnans and Ioffe(1995{\natexlab{b}})}]{bonnans1995second}
Bonnans, Joseph~Fr{\'e}d{\'e}ric, Alexander Ioffe. 1995{\natexlab{b}}.
\newblock Second-order sufficiency and quadratic growth for nonisolated minima.
\newblock {\it Mathematics of Operations Research\/} {\bf 20}(4) 801--817.

\bibitem[{Boyd and Vandenberghe(2009)}]{boyd2009}
Boyd, Stephen, Lieven Vandenberghe. 2009.
\newblock {\it Convex optimization\/}.
\newblock Cambridge university press.

\bibitem[{{Building Robotics}(2016)}]{comfy2016}
{Building Robotics}. 2016.
\newblock Comfy.
\newblock \urlprefix\url{https://gocomfy.com}.

\bibitem[{Burton and Toint(1992)}]{burton1992}
Burton, Didier, Ph~L Toint. 1992.
\newblock On an instance of the inverse shortest paths problem.
\newblock {\it Mathematical Programming\/} {\bf 53}(1-3) 45--61.

\bibitem[{Carr and Lovejoy(2000)}]{carr2000}
Carr, Scott, William Lovejoy. 2000.
\newblock The inverse newsvendor problem: Choosing an optimal demand portfolio
  for capacitated resources.
\newblock {\it Management Science\/} {\bf 46}(7) 912--927.

\bibitem[{Chan et~al.(2014)Chan, Craig, Lee, and Sharpe}]{chan2014}
Chan, Timothy~CY, Tim Craig, Taewoo Lee, Michael~B Sharpe. 2014.
\newblock Generalized inverse multiobjective optimization with application to
  cancer therapy.
\newblock {\it Operations Research\/} .

\bibitem[{Chatterjee(2014)}]{chatterjee2014}
Chatterjee, Sourav. 2014.
\newblock A new perspective on least squares under convex constraint.
\newblock {\it Ann. Statist.\/} {\bf 42}(6) 2340--2381.
\newblock \doi{10.1214/14-AOS1254}.
\newblock \urlprefix\url{http://dx.doi.org/10.1214/14-AOS1254}.

\bibitem[{Crama et~al.(2008)Crama, Reyck, and Degraeve}]{crama2008milestone}
Crama, Pascale, Bert~De Reyck, Zeger Degraeve. 2008.
\newblock Milestone payments or royalties? contract design for r\&d licensing.
\newblock {\it Operations Research\/} {\bf 56}(6) 1539--1552.

\bibitem[{Dempe et~al.(2015)Dempe, Kalashnikov, P{\'e}rez-Vald{\'e}s, and
  Kalashnikova}]{dempe2015}
Dempe, Stephan, Vyacheslav Kalashnikov, Gerardo P{\'e}rez-Vald{\'e}s, Nataliya
  Kalashnikova. 2015.
\newblock {\it Bilevel Programming Problems\/}.
\newblock Springer.

\bibitem[{Efron and Tibshirani(1994)}]{efron1994}
Efron, B., R.J. Tibshirani. 1994.
\newblock {\it An Introduction to the Bootstrap\/}.
\newblock Chapman \& Hall/CRC Monographs on Statistics \& Applied Probability,
  Taylor \& Francis.

\bibitem[{Erkin et~al.(2010)Erkin, Bailey, Maillart, Schaefer, and
  Roberts}]{erkin2010}
Erkin, Zeynep, Matthew~D Bailey, Lisa~M Maillart, Andrew~J Schaefer, Mark~S
  Roberts. 2010.
\newblock Eliciting patients' revealed preferences: An inverse markov decision
  process approach.
\newblock {\it Decision Analysis\/} {\bf 7}(4) 358--365.

\bibitem[{Farag{\'o} et~al.(2003)Farag{\'o}, Szentesi, and
  Szviatovszki}]{farago2003}
Farag{\'o}, Andr{\'a}s, {\'A}ron Szentesi, Bal{\'a}zs Szviatovszki. 2003.
\newblock Inverse optimization in high-speed networks.
\newblock {\it Discrete Applied Mathematics\/} {\bf 129}(1) 83--98.

\bibitem[{Green and Srinivasan(1990)}]{green1990conjoint}
Green, Paul~E, Venkat Srinivasan. 1990.
\newblock Conjoint analysis in marketing: new developments with implications
  for research and practice.
\newblock {\it The Journal of Marketing\/}  3--19.

\bibitem[{Greenshtein and Ritov(2004)}]{greenshtein2004}
Greenshtein, E., Y.~Ritov. 2004.
\newblock Persistence in high-dimensional linear predictor selection and the
  virtue of overparametrization.
\newblock {\it Bernoulli\/} {\bf 10}(6) 971--988.
\newblock \doi{10.3150/bj/1106314846}.
\newblock \urlprefix\url{http://dx.doi.org/10.3150/bj/1106314846}.

\bibitem[{Hastie et~al.(2009)Hastie, Tibshirani, and Friedman}]{hastie2009}
Hastie, T., R.~Tibshirani, J.~Friedman. 2009.
\newblock {\it The Elements of Statistical Learning\/}.
\newblock 2nd ed. Springer-Verlag.

\bibitem[{Haviv and Regev(2012)}]{haviv2012}
Haviv, Ishay, Oded Regev. 2012.
\newblock Tensor-based hardness of the shortest vector problem to within almost
  polynomial factors.
\newblock {\it Theory of Computing\/} {\bf 8}(1) 513--531.

\bibitem[{Heuberger(2004)}]{heuberger2004}
Heuberger, Clemens. 2004.
\newblock Inverse combinatorial optimization: A survey on problems, methods,
  and results.
\newblock {\it Journal of Combinatorial Optimization\/} {\bf 8}(3) 329--361.

\bibitem[{Hillar and Lim(2013)}]{hillar2013}
Hillar, C., L.-H. Lim. 2013.
\newblock Most tensor problems are np-hard.
\newblock {\it J. ACM\/} {\bf 60}(6) 45:1--45:39.
\newblock \doi{10.1145/2512329}.
\newblock \urlprefix\url{http://doi.acm.org/10.1145/2512329}.

\bibitem[{Hochbaum(2003)}]{hochbaum2003}
Hochbaum, Dorit~S. 2003.
\newblock Efficient algorithms for the inverse spanning-tree problem.
\newblock {\it Operations Research\/} {\bf 51}(5) 785--797.

\bibitem[{Iyengar and Kang(2005)}]{iyengar2005}
Iyengar, Garud, Wanmo Kang. 2005.
\newblock Inverse conic programming with applications.
\newblock {\it Operations Research Letters\/} {\bf 33}(3) 319--330.

\bibitem[{Jennrich(1969)}]{jennrich1969}
Jennrich, Robert~I. 1969.
\newblock Asymptotic properties of non-linear least squares estimators.
\newblock {\it The Annals of Mathematical Statistics\/}  633--643.

\bibitem[{Jos\'{e} Fortuny-Amat(1981)}]{fortuny1981}
Jos\'{e} Fortuny-Amat, Bruce~McCarl. 1981.
\newblock A representation and economic interpretation of a two-level
  programming problem.
\newblock {\it The Journal of the Operational Research Society\/} {\bf 32}(9)
  783--792.

\bibitem[{Keshavarz et~al.(2011)Keshavarz, Wang, and Boyd}]{keshavarz2011}
Keshavarz, Arezou, Yang Wang, Stephen Boyd. 2011.
\newblock Imputing a convex objective function.
\newblock {\it Intelligent Control (ISIC), 2011 IEEE International Symposium
  on\/}. IEEE, 613--619.

\bibitem[{Ratliff et~al.(2014{\natexlab{a}})Ratliff, Dong, Ohlsson, and
  Sastry}]{ratliff2014}
Ratliff, Lillian~J., Roy Dong, Henrik Ohlsson, S.~Shankar Sastry.
  2014{\natexlab{a}}.
\newblock Incentive design and utility learning via energy disaggregation.
\newblock {\it 19th World Congress of the International Federation of Automatic
  Control\/}.

\bibitem[{Ratliff et~al.(2014{\natexlab{b}})Ratliff, Dong, Ohlsson, and
  Sastry}]{ratliff2013incentive}
Ratliff, Lillian~J, Roy Dong, Henrik Ohlsson, S~Shankar Sastry.
  2014{\natexlab{b}}.
\newblock Incentive design and utility learning via energy disaggregation.
\newblock {\it Proceedings of the 19th IFAC World Congress\/}. 3158--3163.

\bibitem[{Rockafellar and Wets(1998)}]{rockafellar1998}
Rockafellar, R~Tyrrell, Roger J-B Wets. 1998.
\newblock {\it Variational analysis\/}, vol. 317.
\newblock Springer.

\bibitem[{Saez-Gallego et~al.(2016)Saez-Gallego, Morales, Zugno, and
  Madsen}]{saez2016}
Saez-Gallego, J., J.~M. Morales, M.~Zugno, H.~Madsen. 2016.
\newblock A data-driven bidding model for a cluster of price-responsive
  consumers of electricity.
\newblock {\it IEEE Transactions on Power Systems\/} {\bf PP}(99) 1--11.
\newblock \doi{10.1109/TPWRS.2016.2530843}.

\bibitem[{Schaefer(2009)}]{schaefer2009}
Schaefer, Andrew~J. 2009.
\newblock Inverse integer programming.
\newblock {\it Optimization Letters\/} {\bf 3}(4) 483--489.

\bibitem[{Tao(2012)}]{tao2012}
Tao, T. 2012.
\newblock {\it Topics in Random Matrix Theory\/}.
\newblock Graduate studies in mathematics, American Mathematical Society.

\bibitem[{Troutt et~al.(2006)Troutt, Pang, and Hou}]{troutt2006}
Troutt, Marvin~D, Wan-Kai Pang, Shui-Hung Hou. 2006.
\newblock Behavioral estimation of mathematical programming objective function
  coefficients.
\newblock {\it Management science\/} {\bf 52}(3) 422--434.

\bibitem[{Tversky and Kahneman(1981)}]{tversky1981}
Tversky, Amos, Daniel Kahneman. 1981.
\newblock The framing of decisions and the psychology of choice.
\newblock {\it Science\/} {\bf 211}(4481) 453--458.

\bibitem[{van~der Vaart(2000)}]{van2000}
van~der Vaart, A.W. 2000.
\newblock {\it Asymptotic Statistics\/}.
\newblock Cambridge Series in Statistical and Probabilistic Mathematics,
  Cambridge University Press.

\bibitem[{Vershynin(2012)}]{vershynin2012}
Vershynin, Roman. 2012.
\newblock {\it Compressed Sensing\/}, chap. Introduction to the non-asymptotic
  analysis of random matrices.
\newblock Cambridge University Press, 210--268.

\bibitem[{Wald(1949)}]{wald1949}
Wald, Abraham. 1949.
\newblock Note on the consistency of the maximum likelihood estimate.
\newblock {\it Ann. Math. Statist.\/} {\bf 20}(4) 595--601.
\newblock \doi{10.1214/aoms/1177729952}.
\newblock \urlprefix\url{http://dx.doi.org/10.1214/aoms/1177729952}.

\bibitem[{Wang(2009)}]{wang2009}
Wang, Lizhi. 2009.
\newblock Cutting plane algorithms for the inverse mixed integer linear
  programming problem.
\newblock {\it Operations Research Letters\/} {\bf 37}(2) 114--116.

\bibitem[{Zhang and Zenios(2008)}]{zhang2008dynamic}
Zhang, Hao, Stefanos Zenios. 2008.
\newblock A dynamic principal-agent model with hidden information: Sequential
  optimality through truthful state revelation.
\newblock {\it Operations Research\/} {\bf 56}(3) 681--696.

\bibitem[{Zhang and Liu(1996)}]{zhang1996}
Zhang, Jianzhong, Zhenhong Liu. 1996.
\newblock Calculating some inverse linear programming problems.
\newblock {\it Journal of Computational and Applied Mathematics\/} {\bf 72}(2)
  261--273.

\bibitem[{Zhang and Xu(2010)}]{zhang2010}
Zhang, Jianzhong, Chengxian Xu. 2010.
\newblock Inverse optimization for linearly constrained convex separable
  programming problems.
\newblock {\it European Journal of Operational Research\/} {\bf 200}(3)
  671--679.

\end{thebibliography}

\newpage

\begin{APPENDIX}{} 

\section{Lemmas and Omitted Proofs}

\begin{lemma}
\label{lemma:1}
Suppose $\mathbf{R4}$ holds.  Then for $t > c_1\cdot \gamma$ we have
\begin{equation}
\textstyle\mathbb{P}\Big(\big|\gamma^{-m}\cdot\frac{1}{n}\sum_{j=1}^nK\big(\frac{u_j-u_i}{\gamma}\big) - \mu(u_i)\big|>t\Big) \leq 2\exp\Big(-2c_2 n\gamma^{2m}\cdot(t-c_1\cdot \gamma)^2\Big),
\end{equation}
where $c_1,c_2>0$ are constants.
\end{lemma}

\proof{Proof. }Recall $\mu(u)$ is the probability density function of $u$, and note that
\begin{equation}
\label{eqn:bias1}
\begin{aligned}
\big|\mu(u_i) - \mathbb{E}\big[\gamma^{-m}K\big(\textstyle\frac{u-u_i}{\gamma}\big) \big| u_i\big]\big| &= \textstyle\big|\mu(u_i) - \gamma^{-m}\int_{\mathbb{R}^m}K\big(\textstyle\frac{u-u_i}{\gamma}\big)\mu(u)du\big|\\
&= \textstyle\big|\mu(u_i) - \gamma^{-m}\int_{\mathbb{R}^m} K(s)\mu(u_i+\gamma s)\gamma^mds\big|\\
&= \textstyle\big|\mu(u_i) - \int_{\mathbb{R}^m} K(s)\big(\mu(u_i)+\gamma\nabla \mu(u_i + \beta \gamma s)^Ts\big)ds\big|\\
&= \textstyle\big|\int_{\mathbb{R}^m} K(s)\nabla \mu(u_i + \beta \gamma s)^Tsds\big|\cdot \gamma\\
&\leq c_1 \cdot \gamma,
\end{aligned}
\end{equation}
where the second line follows from a change of variables $s = (u-u_i)/\gamma$, the third line follows from the multivariate form of Taylor's Theorem with some $\beta \in [0,1]$, the fourth line follows because a Kernel function has the property $\int K(u)du = 1$, and the fifth line follows by setting $c_1 = \max_{u\in\mathcal{U}}|\int_{\mathbb{R}^m} K(s)\nabla \mu(u)^Tsds|$.  Note this $c_1$ term is finite because (i) a kernel function has the property that its support is finite (i.e., $K(u) = 0$ for $\|u\| > 1$), and (ii) $\mu(u)$ is a continuously differentiable probability density function by \textbf{R4}.  Next, note that by Hoeffding's inequality \citep{vershynin2012} we have for $t >0$ that
\begin{equation}
\label{eqn:hoeff1}
\textstyle\mathbb{P}\Big(\big|\gamma^{-m}\cdot\frac{1}{n}\sum_{j=1}^nK\big(\frac{u_j-u_i}{\gamma}\big) - \mathbb{E}\big[\gamma^{-m}K\big(\frac{u-u_i}{\gamma}\big)\big|u_i\big]\big|>t\Big) \leq 2\exp\Big(-2c_2 n\gamma^{2m}t^2\Big),
\end{equation}
where $c_2 = (\max_u K(u))^2$.  Combining (\ref{eqn:bias1}) and (\ref{eqn:hoeff1}) gives the desired result.
\Halmos\endproof

\begin{lemma}
\label{lemma:2}
Suppose $\mathbf{A1}$ and $\mathbf{R1}$--$\mathbf{R4}$ hold.  Then for $t > c_3\cdot \gamma^{1/2} + c_4\cdot \gamma$ we have
\begin{equation}
\textstyle\mathbb{P}\Big(\big\|\gamma^{-m}\cdot\frac{1}{n}\sum_{j=1}^ny_j\cdot K\big(\frac{u_j-u_i}{\gamma}\big) - \mu(u_i)\mathcal{S}(u_i,\theta_0)\big\|>t\Big) \leq 2\exp\Big(-2c_5 n\gamma^{2m}\cdot(t-c_3\cdot \gamma^{1/2}-c_4\cdot \gamma)\Big).
\end{equation}
where $c_3,c_4,c_5>0$ are constants.
\end{lemma}

\proof{Proof. }First, note that $\mathcal{S}(u,\theta)$ consists of a single point from the strict convexity assumption in \textbf{R3}. Next, note that having $\mathbf{A1}$ and $\mathbf{R1}$--$\mathbf{R4}$ means that Proposition 4.41 of \citep{bonnans2000} holds: This means for $\gamma > 0$ sufficiently small we have 
\begin{equation}
\label{eqn:lip}
\|\mathcal{S}(u,\theta_0) - \mathcal{S}(u_i,\theta_0)\|\leq\alpha\cdot \gamma^{1/2},
\end{equation}
where $\alpha > 0$ is a constant, whenever $\|u-u_i\|\leq \gamma$.  Next, recall that $y_i$ conditioned on $u_i$ has distribution $\mathcal{S}(u_i,\theta_0)+w_i$ under $\textbf{IC}$.  Moreover, we have
\begin{equation}
\mathbb{E}\big[\gamma^{-m}yK\big(\textstyle\frac{u-u_i}{\gamma}\big) \big| u_i\big] = \mathbb{E}\big[\gamma^{-m}\mathcal{S}(u,\theta_0)K\big(\textstyle\frac{u-u_i}{\gamma}\big) \big| u_i\big],
\end{equation}
since $\mathbb{E}(w_i) = 0$ and $w_i$ is independent of $u_i$.  Thus, we have
\begin{equation}
\label{eqn:bias2}
\begin{aligned}
&\big\|\mu(u_i)\mathcal{S}(u_i,\theta_0) - \mathbb{E}\big[\gamma^{-m}yK\big(\textstyle\frac{u-u_i}{\gamma}\big) \big| u_i\big]\big\| \\
&\quad = \textstyle\big\|\mu(u_i)\mathcal{S}(u_i,\theta_0) - \gamma^{-m}\int_{\mathbb{R}^m}K\big(\textstyle\frac{u-u_i}{\gamma}\big)\mu(u)\mathcal{S}(u,\theta_0)du\big\|\\
&\quad= \textstyle\big\|\mu(u_i)\mathcal{S}(u_i,\theta_0) - \gamma^{-m}\int_{\mathbb{R}^m} K(s)\mu(u_i+\gamma s)\mathcal{S}(u_i+\gamma s,\theta_0)\gamma^mds\big\|\\
&\quad= \textstyle\big\|\mu(u_i)\mathcal{S}(u_i,\theta_0) - \int_{\mathbb{R}^m} K(s)\big(\mu(u_i)+\gamma\nabla \mu(u_i + \beta \gamma s)^Ts\big)\big(\mathcal{S}(u_i,\theta_0) + \\
&\qquad\qquad\qquad\qquad\qquad\qquad\qquad\qquad\qquad\qquad\qquad\qquad\mathcal{S}(u_i+\gamma s,\theta_0) - \mathcal{S}(u_i,\theta_0)\big)ds\big\|\\
&\quad= \textstyle\big\|\int_{\mathbb{R}^m} K(s)\mu(u_i)\big(\mathcal{S}(u_i+\gamma s,\theta_0) - \mathcal{S}(u_i,\theta_0)\big)ds + \int_{\mathbb{R}^m} K(s)\gamma\nabla \mu(u_i + \beta \gamma s)^Ts\mathcal{S}(u,\theta_0)ds\big\|\\
&\quad\leq c_3 \cdot \gamma^{1/2} + c_4\cdot \gamma,
\end{aligned}
\end{equation}
where the second line follows from a change of variables $s = (u-u_i)/\gamma$, the third line follows from the multivariate form of Taylor's Theorem with some $\beta \in [0,1]$, the fourth line follows because a Kernel function has the property $\int K(u)du = 1$, and the fifth line follows from (\ref{eqn:lip}) and by setting $c_3 = \alpha\cdot\max_{u\in\mathcal{U}}|\int_{\mathbb{R}^m} K(s)\mu(u)ds|$ and $c_4 = \max_{u\in\mathcal{U}}(|\int_{\mathbb{R}^m} K(s)\nabla \mu(u)^Tsds|\cdot\|\mathcal{S}(u,\theta_0)\|)$.  Note the $c_3,c_4$ terms are finite because (i) a kernel function has the property that its support is finite (i.e., $K(u) = 0$ for $\|u\| > 1$), (ii) $\mu(u)$ is a continuously differentiable probability density function by \textbf{R4}, and (iii) $\mathcal{S}(u,\theta_0)$ is bounded by \textbf{R1}.  Next, note that $y$ is a sub-exponential random variable \citep{vershynin2012} since (i) $\mathcal{S}(u,\theta_0)$ is a bounded random variable by \textbf{R1}, and (ii) $w$ is sub-exponential by \textbf{R4}.  Hence, by Hoeffding's inequality for sub-exponential random variables \citep{vershynin2012} we have for $t > 0$ that
\begin{equation}
\label{eqn:hoeff2}
\textstyle\mathbb{P}\Big(\big\|\gamma^{-m}\cdot\frac{1}{n}\sum_{j=1}^ny_j\cdot K\big(\frac{u_j-u_i}{\gamma}\big) - \mathbb{E}\big[\gamma^{-m}yK\big(\frac{u-u_i}{\gamma}\big)\big|u_i\big]\big\|>t\Big) \leq 2\exp\Big(-2c_5 n\gamma^{2m}t\Big),
\end{equation}
for some $c_5 > 0$.  Combining (\ref{eqn:bias2}) and (\ref{eqn:hoeff2}) gives the desired result.
\Halmos\endproof 

\vspace{4mm}

\proof{Proof of Proposition 1.} We show this using a counterexample.  Suppose $\mathsf{FOP}$ is $\min\{x^2-(\theta+u)\cdot x\ |\ x\in[0,10]\}$, and note its solution set $\mathcal{S}(u,\theta) = \min\{\frac{u+\theta}{2},10\}$ is single-valued.  Assume the distribution of $u$ is
\begin{equation}
u = \begin{cases} \hphantom{2}0, &\text{with probability (w.p.) } \frac{1}{2}\\ 20, &\text{w.p. } \frac{1}{2}\end{cases}
\end{equation}
and that the distribution of $w$ is
\begin{equation}
w = \begin{cases} -1, &\text{w.p. } \frac{1}{2}\\ +1, &\text{w.p. } \frac{1}{2}\end{cases}
\end{equation}
Finally, suppose $y = \mathcal{S}(u,\theta) + w$, $\Theta = \{\theta\in\mathbb{R} : 0 \leq \theta \leq 10\}$, and $\theta_0 = 10$.  By construction, this problem satisfies $\mathbf{A1}$,$\mathbf{A2}$,$\mathbf{IC}$.  Also, observe that the joint distribution of $(u,y)$ is
\begin{equation}
(u,y) = \begin{cases} (\hphantom{2}0,\hphantom{2}4), &\text{w.p. } \frac{1}{4}\\ (\hphantom{2}0,\hphantom{2}6), &\text{w.p. } \frac{1}{4}\\(20,\hphantom{2}9), &\text{w.p. } \frac{1}{4}\\(20,11), &\text{w.p. } \frac{1}{4}\end{cases}
\end{equation}
We show that both $\mathsf{VIA}$ and $\mathsf{KKA}$ are not estimation consistent for this problem.

We begin with $\mathsf{VIA}$.  This approach solves
\begin{equation}
\label{eqn:via}
\begin{aligned}
\min_{\theta\in\Theta}\ & \textstyle\frac{1}{n}\sum_{i=1}^n\epsilon_i^2\\
\text{s.t. }& \nabla f(y_i,u_i,\theta)\cdot(x_i-y_i) \geq -\epsilon_i, \forall x_i\in[0,10], \qquad\forall i\in[n]\\
\end{aligned}
\end{equation}
The constraint
\begin{equation}
\label{eqn:vi}
\nabla f(y_i,u_i,\theta)\cdot(x_i-y_i) \geq -\epsilon_i, \forall x_i\in[0,10]
\end{equation}
is a variational inequality, and \textsf{VIA} exactly reformulates this using linear duality.  We operate with the original variational inequality since the reformulation in \textsf{VIA} is exact and does not change the solution.  If $y_i = 4$, then a straightforward calculation gives that (\ref{eqn:vi}) is equivalent to the constraint: $\epsilon_i \geq 4\cdot(8-\theta)$ if $\theta \leq 8$, and $\epsilon_i \geq -6\cdot(8-\theta)$ if $\theta > 8$.  If $y_i = 6$, then (\ref{eqn:vi}) is equivalent to the constraint $\epsilon_i \geq 6\cdot(12-\theta)$.  If $y_i = 9$, then (\ref{eqn:vi}) is equivalent to the constraint $\epsilon_i \geq 2 + k$.  Finally, if $y_i = 11$, then (\ref{eqn:vi}) is equivalent to the constraint: $\epsilon_i \geq 11\cdot(2-\theta)$ if $\theta \leq 2$, and $\epsilon_i \geq 2-\theta$ if $\theta > 2$.  Next, we solve the problem $\min\{\epsilon_i^2\ |\ (\ref{eqn:vi})\}$ for each possible value of $y_i$ and $\theta$.  If $y_i = 4$, then the minimum is $16\cdot(8-\theta)^2$ if $\theta\leq 8$, and $36\cdot(8-\theta)^2$ if $\theta > 8$.  If $y_i = 6$, then the minimum is $36\cdot(12-\theta)^2$.  If $y_i = 9$, then minimum is $(2+\theta)^2$.  If $y_i = 11$, then the minimum is $121\cdot(2-\theta)^2$ if $\theta \leq 2$, and $0$ if $\theta > 2$.  Thus, we have
\begin{equation}
4\cdot\mathbb{E}(\epsilon_i^2) = \begin{cases}
36\cdot(12-\theta)^2 + (2+\theta)^2+121\cdot(2-\theta)^2+16\cdot(8-\theta)^2, &\text{if } \theta \leq 2\\
36\cdot(12-\theta)^2 + (2+\theta)^2+16\cdot(8-\theta)^2, &\text{if } \theta \in (2,8]\\
36\cdot(12-\theta)^2 + (2+\theta)^2+36\cdot(8-\theta)^2, &\text{if } \theta > 8
\end{cases}
\end{equation}
Finally, we solve the optimization problem $\min\{\mathbb{E}(\epsilon_i^2)\ |\ \theta\in[0,10]\}$.  A simple calculation gives that the minimum occurs at $\theta^* = \frac{718}{73} \approx 9.8356$.  However, the minimizer of (\ref{eqn:via}) will converge in probability to $\theta^*$, because (i) we can exactly reformulate (\ref{eqn:via}) as
\begin{equation}
\label{eqn:via-ref}
\begin{aligned}
\min_{\theta\in\Theta}\ & \textstyle\frac{1}{n}\sum_{i=1}^n\epsilon_i^2\\
\text{s.t. }& \epsilon_i^2 = \begin{cases} 16\cdot(8-\theta)^2\cdot\mathds{1}(\theta \leq 8) + 36\cdot(8-\theta)^2\cdot\mathds{1}(\theta > 8), & \text{if } y_i = 4\\
36\cdot(12-\theta)^2,  & \text{if } y_i = 6\\
(2+\theta)^2,  & \text{if } y_i = 9\\
121\cdot(2-\theta)^2\cdot\mathds{1}(\theta \leq 2),  & \text{if } y_i = 11
\end{cases} \qquad\forall i\in[n]\\
\end{aligned}
\end{equation}
which (ii) implies we can apply the uniform law of large numbers \citep{jennrich1969} since $\epsilon_i^2$ as defined in (\ref{eqn:via-ref}) is a continuous function, and thus (iii) we get convergence of the minimizer from a standard consistency result in statistics (see for instance Theorem 5.7 in \citep{van2000} or Theorem 5.2.3 in \citep{bickel2006}).  This shows \textsf{VIA} is not estimation consistent, since $\theta_0 = 10$.

Next, we consider \textsf{KKA}.  This approach solves
\begin{equation}
\label{eqn:kka}
\begin{aligned}
\min_{\theta\in\Theta}\ & \textstyle\frac{1}{n}\sum_{i=1}^n\|\epsilon_i\|^2\\
\text{s.t. }& \nabla f(y_i,u_i,\theta) -  \lambda_{i1} + \lambda_{i2} =  \epsilon_{i1}\\
&- \lambda_{i1}\cdot y_i =  \epsilon_{i2}\\
& \lambda_{i2}\cdot (y_i -10)=  \epsilon_{i3}\\
&\lambda_i \geq 0\\
\end{aligned}
\end{equation}
We first solve the problem (\ref{eqn:kka}), with $n=1$, for each possible value of $y_i$ and $\theta$.  If $y_i = 4$, then the minimum is $\frac{16}{17}\cdot(8-\theta)^2$ if $\theta \leq 8$, and $\frac{36}{37}\cdot(8-\theta)^2$ if $\theta > 8$.  If $y_i = 6$, then the minimum is $\frac{36}{37}\cdot(12-\theta)^2$.  If $y_i = 9$, then the minimum is $\frac{1}{2}\cdot(2+\theta)^2$.  If $y_i = 11$, then the minimum is $\frac{121}{122}\cdot(2-\theta)^2$ if $\theta \leq 2$, and $\frac{1}{2}\cdot(2-\theta)^2$ if $\theta > 2$.  Thus, we have
\begin{equation}
4\cdot\mathbb{E}(\|\epsilon_i\|^2) = \begin{cases}
\frac{36}{37}\cdot(12-\theta)^2 + \frac{1}{2}\cdot(2+\theta)^2+\frac{121}{122}\cdot(2-\theta)^2+\frac{16}{17}\cdot(8-\theta)^2, &\text{if } \theta \leq 2\\
\frac{36}{37}\cdot(12-\theta)^2 + \frac{1}{2}\cdot(2+\theta)^2+\frac{1}{1}\cdot(2-\theta)^2+\frac{16}{17}\cdot(8-\theta)^2, &\text{if } \theta \in (2,8]\\
\frac{36}{37}\cdot(12-\theta)^2 + \frac{1}{2}\cdot(2+\theta)^2+\frac{1}{2}\cdot(2-\theta)^2+\frac{36}{37}\cdot(8-\theta)^2, &\text{if } \theta > 8
\end{cases}
\end{equation}
Finally, we solve the optimization problem $\min\{\mathbb{E}(\|\epsilon_i\|^2)\ |\ \theta\in[0,10]\}$.  A simple calculation gives that the minimum occurs at $\theta^* = \frac{12080}{1833} \approx 6.5903$.  However, the minimizer of (\ref{eqn:kka}) will converge in probability to $\theta^*$, because (i) we can exactly reformulate (\ref{eqn:kka}) as
\begin{equation}
\label{eqn:kka-ref}
\begin{aligned}
\min_{\theta\in\Theta}\ & \textstyle\frac{1}{n}\sum_{i=1}^n\|\epsilon_i\|^2\\
\text{s.t. }& \|\epsilon_i\|^2 = \begin{cases} \frac{16}{17}\cdot(8-\theta)^2\cdot\mathds{1}(\theta \leq 8) + \frac{36}{37}\cdot(8-\theta)^2\cdot\mathds{1}(\theta > 8), & \text{if } y_i = 4\\
\frac{36}{37}\cdot(12-\theta)^2,  & \text{if } y_i = 6\\
\frac{1}{2}\cdot(2+\theta)^2,  & \text{if } y_i = 9\\
\frac{121}{122}\cdot(2-\theta)^2\cdot\mathds{1}(\theta \leq 2) + \frac{1}{2}\cdot(2-\theta)^2,  & \text{if } y_i = 11
\end{cases} \qquad\forall i\in[n]\\
\end{aligned}
\end{equation}
which (ii) implies we can apply the uniform law of large numbers \citep{jennrich1969} since $\|\epsilon_i\|^2$ as defined in (\ref{eqn:kka-ref}) is a continuous function, and thus (iii) we get convergence of the minimizer from a standard consistency result in statistics (see for instance Theorem 5.7 in \citep{van2000} or Theorem 5.2.3 in \citep{bickel2006}).  This shows that \textsf{KKA} is not estimation consistent, since $\theta_0 = 10$.
\Halmos\endproof

\vspace{4mm}

\proof{Proof of Corollary 2.} It suffices to show that risk consistency is necessary for estimation consistency in the counterexample given in the proof of Proposition \ref{proposition:estincon}. First note that the risk function
\begin{equation}
Q(\theta) = \textstyle\mathbb{E}\Big(\|y-\min\{\frac{u+\theta}{2},10\}\|^2\Big) = \frac{1}{4}\cdot\Big((4-\frac{\theta}{2})^2 + (6-\frac{\theta}{2})^2+ (9-10)^2 +(11-10)^2\Big)
\end{equation}
is continuous since $\Theta = \{\theta\in\mathbb{R} : 0 \leq \theta\leq 10\}$.  Now suppose a sequence $\hat{\theta}_n$ is estimation consistent. Since $\hat{\theta}_n \stackrel{p}{\longrightarrow} \theta_0$, by continuity of $Q(\theta)$ and the continuous mapping theorem \citep{van2000}, we have $Q(\hat{\theta}_n) \stackrel{p}{\longrightarrow} Q(\theta_0)$. Since $\arg\min\{Q(\theta)\ |\ \theta\in\Theta\} = 10 = \theta_0$, and $\hat{\theta}_n \longrightarrow \theta_0$ we have that $\hat{\theta}_n$ converges to a minimizer of $Q(\theta)$. Hence $\hat{\theta}_n$ is risk consistent. \Halmos\endproof

\vspace{4mm}

\section{Identifiability in Inverse Optimization}
\label{sect:iden}
Estimation consistency in any statistical setting (including inverse optimization with noisy data) requires that an identifiability condition holds, and such identifiability conditions can be stated under a variety of different mathematical formulations \citep{wald1949,jennrich1969,bartlett2002,greenshtein2004,bickel2006,chatterjee2014,aswani2015}.  The intuition for these different formulations is the same: Essentially, an identifiability condition states that the output of the model is different for two distinct sets of model parameters. It is important to note that identifiability is a statistical property of the model and the error metric used.  Consequently, it is possible for an estimator to be statistically inconsistent, even when an identifiability condition holds (see for instance Proposition {proposition:estincon}).  In the context of inverse optimization with noisy data, we define an identifiability condition \textbf{IC}. 

 Showing that \textbf{IC} holds is complicated by the presence of constraints in \textsf{FOP}. To illustrate this, consider two related instances of \textsf{FOP} with $x\in\mathbb{R}$ and $\theta \in[0,2]$.  The first $\min (x-\theta)^2$ is \textsf{FOP-I}, and the second $\min \{(x-\theta)^2\ |\ x\leq 1\}$ is \textsf{FOP-II}.  Since these two problems are strictly convex, their minimizers are unique.  Next, suppose we would like to estimate $\theta$ given a (noiseless) measurement $y_i$ of the minimizer.  Observe that \textsf{FOP-I} is identifiable because we must have $\theta = y_i$.  However, \textsf{FOP-II} is not identifiable because if $y_i = 1$, then we may have any $\theta \in [1,2]$.  Thus, the constraint $x \leq 1$ renders \textsf{FOP-II} unidentifiable, and precludes the possibility of \textbf{IC} holding for \textsf{FOP-II}.

Though \textsf{FOP-II} is not identifiable, a related problem is identifiable because of external inputs.  In particular, consider an \textsf{FOP-III}  with $x\in\mathbb{R}$ and $\theta \in[0,2]$ that is given by $\min \{(x - \theta - u)^2\ |\ x \leq 1\}$.  This problem is strictly convex, and so its minimizer is unique for each fixed value of $u$.  In fact, the minimizer is given by $y_i = \min\{(\theta+u_i),1\}$.  And so a sufficient condition for identifiability of \textsf{FOP-III} is if $\mathbb{P}(u_i \leq -1) > 0$.  For instance, if $u_i = -1$ then $y_i = \theta - 1$ and so $\theta$ is uniquely determined by $y_i$.  The presence of the input parameter $u$ ensures identifiability of \textsf{FOP-III}.

\end{APPENDIX}



\end{document}